\documentclass[a4paper, 12pt, twoside]{article}

\usepackage[a4paper, margin=1in]{geometry}
\usepackage[latin1]{inputenc}
\usepackage[T1]{fontenc}
\usepackage[english]{babel}
\usepackage{textcomp} 
\usepackage{graphicx}
\DeclareGraphicsExtensions{.jpg,.mps,.pdf,.png,.gif}
\usepackage[french]{varioref} 
\usepackage[]{amsmath,amssymb,amsfonts}

\usepackage{float}
\usepackage{bbm}
\usepackage{bbold}
\usepackage{amsthm}
\usepackage{dsfont}
\usepackage{color}
\usepackage{natbib}
\usepackage{pict2e}
\usepackage{color, colortbl}

\definecolor{blue1}{rgb}{0,0,0.5451} 
  \definecolor{blue2}{rgb}{0.2,0.4,0.65}
\definecolor{DarkBlue}{rgb}{0,0,0.5451}  
\definecolor{darkred}{rgb}{0.75,0.0,0.0}
\definecolor{darkgreen}{rgb}{0.0,0.6,0.0}
\definecolor{darkblue}{rgb}{0.0,0.0,0.6}
\definecolor{darkcyan}{rgb}{0.0,0.6,0.6}
\definecolor{darkmagenta}{rgb}{0.6,0.0,0.6}
\definecolor{darkyellow}{rgb}{0.6,0.6,0.0}
\definecolor{lightred}{rgb}{1.0,0.9,0.9}
\definecolor{lightgreen}{rgb}{0.9,1.0,0.9}
\definecolor{lightblue}{rgb}{0.9,0.9,1.0}
\definecolor{lightcyan}{rgb}{0.8,1.0,1.0}
\definecolor{lightmagenta}{rgb}{1.0,0.8,1.0}
\definecolor{lightyellow}{rgb}{1.0,1.0,0.8}
\definecolor{paleyellow}{rgb}{1.00,1.0,0.80}
\definecolor{amber}{rgb}{1.0,0.8,0.0}
\definecolor{darkamber}{rgb}{1.0,0.5,0.0} 
\definecolor{Gray}{gray}{0.9}

\newcommand{\Ben}{\begin{enumerate}}
\newcommand{\Een}{\end{enumerate}}
\newcommand{\Bit}{\begin{itemize}}
\newcommand{\Eit}{\end{itemize}}
\newcommand{\Beq}{\begin{equation}}
\newcommand{\Eeq}{\end{equation}}
\newcommand{\Ba}{\begin{align*}}
\newcommand{\Ea}{\end{align*}}
\newcommand{\Mb}{\mathbf}

\newtheorem{Th}{Theorem}

\newtheorem{Prop}{Proposition}

\newtheorem{Corr}{Corollary}
\newtheorem{Def}{Definition}

\title{Space-time max-stable models with spectral separability}
\date{July 17, 2015}
\begin{document}
\author{Paul Embrechts\footnote{ETH Zurich (Department of Mathematics, RiskLab). embrechts@math.ethz.ch}
~~ Erwan Koch\footnote{ETH Zurich (Department of Mathematics, RiskLab). erwan.koch@math.ethz.ch}
~~ Christian Robert\footnote{ISFA Universit\'e Lyon 1. christian.robert@univ-lyon1.fr}
}

\maketitle
\begin{abstract}
Natural disasters may have considerable impact on society as well as on
(re)insurance industry. Max-stable processes are ideally suited for the
modeling of the spatial extent of such extreme events, but it is often
assumed that there is no temporal dependence. Only a few papers have
introduced spatio-temporal max-stable models, extending the Smith, Schlather
and Brown-Resnick spatial processes. These models suffer from two major
drawbacks: time plays a similar role as space and the temporal dynamics is
not explicit. In order to overcome these defects, we introduce
spatio-temporal max-stable models where we partly decouple the influence of
time and space in their spectral representations. We introduce both
continuous and discrete-time versions. We then consider particular Markovian
cases with a max-autoregressive representation and discuss their properties.
Finally, we briefly propose an inference methodology which is tested through
a simulation study.

\medskip

\end{abstract}

\section{Introduction}

In the context of climate change, some extreme events tend to be more and
more frequent; see e.g. \cite{SwissRe}. Meteorological and more generally
environmental disasters have a considerable impact on society as well as the
(re)insurance industry. Hence, the statistical modeling of extremes
constitutes a crucial challenge. Extreme value theory (EVT) provides
powerful statistical tools for this purpose.

EVT can basically be divided into three different streams closely linked to
each other: the univariate case, the multivariate case and the theory of
max-stable processes. For an introduction to the univariate theory, see e.g. 
\cite{coles2001introduction} and for a detailed description, see e.g. \cite%
{embrechts1997modelling} or \cite{beirlant2006statistics}. In the
multivariate case, we refer to \cite{resnickextreme}, \cite%
{beirlant2006statistics} and \cite{de2007extreme}. Max-stable processes
constitute an extension of EVT to the level of stochastic processes %
\citep{haan1984spectral, de1986stationary} and are very well suited for the
modeling of spatial extremes. Indeed, under fairly general conditions, it
can be shown that the distribution of the random field of the suitably normalized temporal maxima at each point of the space is necessarily
max-stable when the number of temporal observations tends to infinity. For a
detailed overview of max-stable processes, we refer to \cite{de2007extreme}.

In the literature about max-stable processes, measurements are often assumed to be independent in time
and thus only the spatial structure is studied 
\citep[see
e.g.][]{padoan2010likelihood}. Nevertheless, the temporal dimension should
be taken into account in a proper way. 
To the best of our knowledge, only a few papers focus on such a question.
The majority of the spatio-temporal models introduced is based on
Schlather's spectral representation \citep{penrose1992semi, schlather2002models} which has given rise to the well-known
Schlather \citep{schlather2002models} and Brown-Resnick %
\citep{kabluchko2009stationary} processes. This representation tells us that if $%
(U_{i})_{i\geq 1}$ generates a Poisson point process on $(0,\infty )$ with
intensity $u^{-2}du$ and $( ( Y_{i}(\mathbf{y) ) }_{\mathbf{y}\in 
\mathds{R}^{d}})_{i\geq 1}$\ are independent and identically
distributed (iid) non-negative stationary stochastic processes such that $\mathds{E}[Y_{i}(%
\mathbf{y)}]=1$ for each $\mathbf{y}\in \mathds{R}^{d}$, then the process $%
(\bigvee_{i=1}^{\infty} \left\{ U_{i} Y_{i}(\mathbf{y)}\right\} )_{\mathbf{y}\in 
\mathds{R}^{d}}$ is stationary simple max-stable, where simple means that the margins are standard Fr\'{e}%
chet. Here $\bigvee$ denotes the max-operator. In \cite{davis2013max}, \cite{huser2014space} and \cite%
{buhl2015anisotropic}, the idea underlying the construction of the
spatio-temporal model is to divide the dimension $d$ into the dimension $d-1$
for the spatial component and the dimension $1$ for the time. \cite%
{davis2013max} introduce the Brown-Resnick model in space and time by taking
a log-normal process for $Y_{i}$ while \cite{buhl2015anisotropic} introduce
an extension of this model to the anisotropic setting. \cite{huser2014space}
consider an extension of the Schlather model by using a truncated Gaussian
process for the $Y_{i}$ and a random set that allows the process to be
mixing in space as well as to exhibit a spatial propagation. Advantages of
these models lie in the facts that the Schlather and Brown-Resnick models
have been widely studied and that the large literature about spatio-temporal
correlation functions for Gaussian processes can be used, allowing for a
considerable diversity of spatio-temporal behavior. \cite{davis2013max} also introduce the spatio-temporal version of the Smith
model \citep{smith1990max} that is based on de Haan's spectral representation \citep[see][]{haan1984spectral}. If $(U_{i},C_{i})_{i\geq 1}$ are the points of a
Poisson point process on $(0,\infty )\times \mathds{R}^{d}$ with intensity $%
u^{-2}du\times dc$ and if $g_{\mathbf{y}}$ are measurable non-negative
functions satisfying $\int_{\mathds{R}^{d}}g_{\mathbf{y}}(c)dc=1$ for each $%
\mathbf{y}\in \mathds{R}^{d}$, then the process $(\bigvee_{i=1}^{\infty} \left\{
U_{i} g_{\mathbf{y}}\left( C_{i}\right) \right\} )_{\mathbf{y}\in \mathds{R}%
^{d}}$ is a simple max-stable process. However, they do not allow any
interaction between the spatial components and the temporal one in the
underlying covariance matrix. The previous spatio-temporal max-stable models suffer from some defects.
First, they are all continuous-time processes whereas measurements in
environmental science are often time-discrete. Second, time has no
specific role but is equivalent to an additional spatial dimension.
Especially, the spatial and temporal distributions belong to a similar class of models.
This constitutes a serious drawback since such a similarity is not supported by any physical argument. Third,
the temporal dynamics is not explicit and hence difficult to identify and
interpret. Finally, these models have in general no causal representation.

The theory of linear ARMA processes has led to the max-autoregressive moving
average processes (MARMA$(p,q)$) introduced by \cite{davis1989basic}. The
real-valued process $(X(t))_{t\in \mathds{Z}}$ follows the MARMA$(p,q)$ model if it
satisfies the recursion 
\begin{equation*}
X(t)=\max (\phi _{1}X(t-1),\dots ,\phi _{p}X(t-p),Z(t),\theta
_{1}Z(t-1),\dots ,\theta _{q}Z(t-q)), \quad t \in \mathds{Z},
\end{equation*}%
where $\phi _{i},\theta _{j}\geq 0$ for $%
i=1,\dots ,p$ and $j=1,\dots ,q$ and the max-stable random variables $Z(t)$
for $t\in \mathds{Z}$ are iid. It is a time series model which is max-stable
in time. However, the spatial aspect is absent. An interesting approach to
build spatio-temporal max-stable processes could be inspired by the theory
of linear processes in function spaces such as Hilbert and Banach spaces %
\citep[see e.g.][]{bosq2000linear} and especially by the autoregressive
Hilbertian model of order 1, ARH(1) 
\citep[see e.g.][]{bosq2000linear,
hormann2012functional}. We say that a sequence $( X(t) )_{t\in \mathds{Z}}$
of mean zero functions in an Hilbert space $H$ follows an ARH(1) process, if 
\begin{equation*}
X(t)=\Psi (X(t-1))+Z(t),\quad t\in \mathds{Z},
\end{equation*}%
where $\Psi $ is a bounded linear operator from $H$ to $H$ and $(Z(t))_{t\in 
\mathds{Z}}$ is a sequence of iid mean zero functions in $H$ satisfying $%
\mathds{E}(\Vert Z(t)\Vert ^{2})<\infty $, where $\Vert .\Vert $ denotes the
norm induced by the scalar product on $H$. Various types of linear
transformations can be applied though the most commonly used is the local
average operator which involves a kernel. A transposition of this model to
the context of the maximum instead of the sum could for instance be written
as 
\begin{equation}
X(t,\mathbf{x})=\max \left( \Psi (X(t-1,\mathbf{\cdot }))(\mathbf{x}),Z(t,%
\mathbf{x})\right) ,\quad (t,\mathbf{x)}\in (\mathds{Z},\mathds{R}^{d-1}),
\label{Eq_Extension_ARH1}
\end{equation}%
where $( ( Z(t,\mathbf{x}) )_{\mathbf{x}\in \mathds{R}^{d-1}})_{t\in \mathds{Z%
}}$ is a sequence of iid spatial max-stable processes and $\Psi $ is an
operator from the space of continuous functions on $\mathds{R}^{d-1}$ to
itself such that, if $(X(t-1,\mathbf{x}))_{\mathbf{x}\in \mathds{R}^{d-1}}$
is max-stable in space, then $\left( \Psi (X(t-1,\mathbf{\cdot }))(\mathbf{x}%
)\right) _{\mathbf{x}\in \mathds{R}^{d-1}}$ is also max-stable in space.
Such an operator could for instance be a \textquotedblleft
moving-maxima\textquotedblright\ operator 
\begin{equation*}
\Psi (X(t-1,\mathbf{\cdot }))(\mathbf{x})=\bigvee_{\mathbf{s}\in \mathds{R}%
^{d-1}}\{K(\mathbf{s},\mathbf{x})X(t-1,\mathbf{s})\},\quad \mathbf{x}\in 
\mathds{R}^{d-1},
\end{equation*}%
where $K$ is a kernel (see \cite{meinguet2012maxima} for a similar idea), or an operator combining a translation in space with a scaling
transformation 
\begin{equation*}
\Psi (X(t-1,\mathbf{\cdot }))(\mathbf{x})=aX(t-1,\mathbf{x}-\boldsymbol{\tau 
}),\quad \mathbf{x}\in \mathds{R}^{d-1},
\end{equation*}%
where $a \in (0,1)$ and $\boldsymbol{\tau }\in \mathds{R}^{d-1}$.

In this paper, we propose a class of models where we partly decouple the
influence of time and space, but such that time influences space through a
bijective operator on space. We present both continuous-time and
discrete-time versions. A first advantage of this class of models lies in
their flexibility since they allow the marginal distribution in time to belong to a different class than the stationary distribution in space. Actually,
these margins can be chosen in function of the application. Due to the spatial operator
mentioned above, our models are able to account for physical processes such
as propagations/contagions/diffusions. Furthermore, the estimation procedure
can be simplified since the purely spatial parameters can be estimated
independently of the purely temporal ones.

Then, we study some particular sub-classes of our general class of models, where the function
related to time in the spectral representation is the exponential density
(in the continuous-time case) or takes as values the probabilities of a
geometric random variable (in the discrete-time case). In this context, our
models become Markovian and have a max-autoregressive representation. This
makes the dynamics of these models explicit and easy to interpret physically.

The remaining of the paper is organized as follows. Section 2 presents our
class of spectrally separable space-time max-stable models. In Section 3, we focus on the particular Markovian cases where the space is $\mathds{R}%
^{2}$ and the unit sphere in $\mathds{R}^{3}$, respectively. Section 4
briefly presents an estimation procedure as well as an application of the
latter on simulated data. Some concluding remarks are given in Section 5.

\section{A new class of space-time max-stable models}

The time index $t$ and space index $\mathbf{x}$ will belong respectively to
the sets $\mathcal{I}$ and $\mathcal{X}$. The models we introduce will be
either continuous-time ($\mathcal{I}=\mathds{R}$) or discrete-time ($%
\mathcal{I}=\mathds{Z}$). In the following, we denote by $\delta $ the
Lebesgue measure on $\mathds{R}$ in the case $\mathcal{I}=\mathds{R}$ and
the counting measure $\sum_{z\in \mathds{Z}}\partial _{\{z\}}$ when $%
\mathcal{I}=\mathds{Z}$, where $\partial $ stands for the Dirac measure.

To define discrete-time models, we need to introduce the notion of
homogeneous Poisson point process on $\mathds{Z}$. Let $\left (
N_{k}\right ) _{k\in \mathds{Z}}$ be iid Poisson$(1)$. For $A\subset \mathds{%
Z}$, $N\left( A\right) =\sum_{k\in A}N_{k}$ defines an homogeneous Poisson
point process on $\mathds{Z}$ with constant intensity equal to one (see
Appendix \ref{Sec_Poisson_Z}). Note that $N$ is not a simple
point process.

Space-time simple max-stable processes on $\mathcal{I}\times \mathcal{X}$ allow for a spectral representation of the
following form (see e.g. \cite{haan1984spectral}):
\begin{equation}
X(t,\mathbf{x})=\bigvee_{i=1}^{\infty } \{ U_{i} V_{(t,\mathbf{x)}}(W_{i}) \}, \quad
(t,\mathbf{x)\in }\mathcal{I}\times \mathcal{X},
\label{GeneralSpectralRepresentation}
\end{equation}%
where $(U_{i},W_{i})_{i\geq 1}$ are the points of a Poisson point process on $(0,\infty
)\times E$ with intensity $u^{-2}du\times \mu (dw)$\ for some Polish measure
space $(E,\mathcal{E},\mu )$ and the functions $V_{(t,\mathbf{x}%
)}:E\rightarrow (0,\infty )$ are measurable such that $\int_{E}V_{(t,%
\mathbf{x})}(w)\mu (dw)=1$ for each $(t,\mathbf{x)\in }\mathcal{I}\times 
\mathcal{X}$. A class of space-time max-stable models avoiding the previously mentioned shortcomings is introduced below.

\begin{Def}[Space-time max-stable models with spectral separability]
\label{Class_Space_Time_Maxstable_Process_Spectral_Separability}
The class of space-time max-stable models with spectral separability is defined inserting the following spectral decomposition in \eqref{GeneralSpectralRepresentation}:
\begin{equation}
V_{(t,\mathbf{x})}(W_{i})=V_{t}(B_{i})V_{R_{(t,B_{i})}\mathbf{x}}(C_{i}),
\label{SpectralRepresentation}
\end{equation}
where:
\begin{itemize}
\item $(U_{i},B_{i},C_{i})_{i\geq 1}$ are the points of a Poisson point process on $%
(0,\infty )\times E_{1}\times E_{2}$ with intensity $u^{-2}du\times \mu
_{1}(db)\times \mu _{2}(dc)$\ for some Polish measure spaces $(E_{1},%
\mathcal{E}_{1},\mu _{1})$ and $(E_{2},\mathcal{E}_{2},\mu _{2})$;

\item the operators $R_{(t,b)}$ are bijective from $\mathcal{X}$ to $%
\mathcal{X}$ for each $(t,b) \in \mathcal{I} \times E_1$;

\item the functions $V_{t}:E_{1}\rightarrow (0,\infty )$ are measurable
such that $\int_{E_{1}}V_{t}(b)\mu _{1}(db)=1$ for each $t\in \mathcal{I}$ and
the functions $V_{\mathbf{x}}:E_{2}\rightarrow (0,\infty )$ are measurable
such that $\int_{E_{2}}V_{\mathbf{x}}(c)\mu _{2}(dc)=1$ for each $\mathbf{x}%
\in \mathcal{X}$.
\end{itemize}
\end{Def}

We emphasize that the models belonging to this class are max-stable in space and time, since 
\begin{align*}
\int_{E}V_{(t,\mathbf{x})}(w)\mu (dw) &=\int_{E_{1}\times
E_{2}}V_{t}(b)V_{R_{(t,b)}\mathbf{x}}(c)\mu _{1}(db)\mu _{2}(dc) 
\\& =\int_{E_{1}}V_{t}(b)\left( \int_{E_{2}}V_{R_{(t,b)}\mathbf{x}}(c)\mu
_{2}(dc)\right) \mu _{1}(db) 
\\&=\int_{E_{1}}V_{t}(b)\mu _{1}(db)=1,
\end{align*}
but of course also in space and in time only. A spectral decomposition in
space e.g. is easily derived since, for a fixed $t$, $%
(U_{i}V_{t}(B_{i}),C_{i})_{i\geq 1}$ defines a Poisson point process on $%
(0, \infty )\times E_{2}$ with intensity $u^{-2}du\times \mu _{2}(dc)$\ and $%
\int_{E_{2}}V_{R_{(t,b)}\mathbf{x}}(c)\mu _{2}(dc)=1$ for each $\mathbf{x\in }%
\mathcal{X}$ and $b\in E_1$.

The crucial point in the previous definition lies in the fact that we have
decoupled the spectral functions with respect to time and the spectral
functions with respect to space given time. This allows one to deal with the
temporal and the spatial aspects separately. Moreover, the latter depends on
time through a bijective transformation which typically may account for an
underlying physical process.

The finite dimensional distributions of $X$ in \eqref{SpectralRepresentation} are given, for $M \in \mathds{N} \backslash \{ 0 \}$, $t_{1},\dots ,t_{M} \in \mathcal{I}$, $ 
\mathbf{x}_{1},\dots ,\mathbf{x}_{M} \in \mathcal{X}$ and $
z_{1},\dots ,z_{M} >0$, by
\begin{align}
& -\log \left( \mathds{P}(X(t_{1},\mathbf{x}_{1})\leq z_{1},\dots ,X(t_{M},%
\mathbf{x}_{M})\leq z_{M})\right) \nonumber
\\& =\int_{E_{1}\times E_{2}}\bigvee_{m=1}^{M}%
\left \{ \frac{V_{t_{m}}(b)V_{R_{(t_{m},b)}\mathbf{x}_{m}}(c)}{z_{m}} \right \} \mu_{1}(db)\mu
_{2}(dc).
\label{Eq_Finite_Dimensional_Distribution}
\end{align}

We now provide some examples of sub-classes of the general class of space-time max-stable processes given in Definition \ref{Class_Space_Time_Maxstable_Process_Spectral_Separability}.

\medskip

\noindent \textbf{i) Models of type 1: de Haan's representation with }$%
\mathcal{X}=\mathds{R}^{2}$

\medskip

We take $E_{1}=\mathcal{I}$ with $\mu
_{1}=\delta $ and $E_{2}=\mathcal{X}=\mathds{R}^{2}$ with $\mu _{2}=\lambda _{2}
$, where $\lambda _{2}$ is the Lebesgue measure on $\mathds{R}^{2}$. Let $g$ be a probability
density function (case $\mathcal{I}=\mathds{R}$) or a discrete probability
distribution (case $\mathcal{I}=\mathds{Z}$), and $f$ be a probability density
function on $\mathds{R}^{2}$. We then assume that%
\begin{equation*}
V_{t}(b)=g(t-b), \quad \quad V_{\mathbf{x}}(c)=f(\mathbf{x}-c)
\end{equation*}%
and that the operators $R_{(t,b)}$ are translations: for all $t,b \in \mathcal{I}$ and $\Mb{x} \in \mathds{R}^2$, $R_{(t,b)}\mathbf{x}=\mathbf{x}%
-(t-b)\boldsymbol{\tau }$, where $\boldsymbol{\tau} \in \mathds{R}^{2}$.

The class of moving maxima max-stable processes with general spectral
representation $\left( \ref{GeneralSpectralRepresentation}\right)$ assumes
the existence of a probability density function $h$ on $\mathcal{I}\times 
\mathds{R}^{2}$ such that 
\begin{equation*}
V_{(t,\mathbf{x})}\left( w\right) =h\left( t-b,\mathbf{x}-c\right). 
\end{equation*}%
 The density function $h$ can always be decomposed as follows: 
\begin{equation*}
h\left( t,\mathbf{x}\right) =g\left( t\right) h_1 \left( \mathbf{x}|t\right),
\end{equation*}%
where $h_1 \left( \mathbf{x}|t\right)$ is the conditional probability density function on $\mathds{R}^{2}$
given $t$. For models of type 1, we have implicitly assumed that this
density function satisfies the equality $h_1 \left( \mathbf{x}|t\right)
=f\left( \mathbf{x}-t\boldsymbol{\tau }\right) $.

Models of this type are interesting in practice since, as we will see in the next section, they have a max-autoregressive representation for a well chosen function $g$. The latter makes the dynamics explicit. Moreover, the translation operator allows to model physical processes such as propagation and diffusion.

We denote by $\mathds{S}^{2}= \{ \Mb{x} \in \mathds{R}^3: \| \Mb{x} \|=1 \}$, the unit sphere in $\mathds{R}^3$.

\medskip

\noindent \textbf{ii) Models of type 2: de Haan's representation with }$%
\mathcal{X}=\mathds{S}^{2}$

\medskip

We choose $E_{1}=\mathcal{I}$ with $\mu
_{1}=\delta $ and $E_{2}=\mathcal{X}=\mathds{S}^{2}$ with $\mu _{2}=\lambda _{%
\mathds{S}^2}$, where $\lambda _{%
\mathds{S}^2}$ is the Lebesgue measure on $\mathds{S}^{2}$. Let $g$ be a
probability density function (case $\mathcal{I}=\mathds{R}$) or a discrete
probability distribution (case $\mathcal{I}=\mathds{Z}$) and $f$ be the von
Mises--Fisher probability density function on $\mathds{S}^{2}$ with
parameters $\boldsymbol{\mu }\in \mathds{S}^{2}$ and $\kappa \geq 0$:
\begin{equation}
\label{Eq_von_Mises_Fisher}
f(\mathbf{x};\boldsymbol{\mu },\kappa )=\frac{\kappa }{4\pi \sinh \kappa }\exp
\left( \kappa \boldsymbol{\mu }' \mathbf{x}\right) \text{,}\quad 
\mathbf{x}\in \mathds{S}^{2}\text{.}
\end{equation}%
The parameters $\boldsymbol{\mu }$ and $\kappa $ are called the mean direction
and concentration parameter, respectively. The greater the value of $\kappa $%
, the higher the concentration of the distribution around the mean direction 
$\boldsymbol{\mu }$. The distribution is uniform on the sphere for $\kappa =0$ and unimodal for $\kappa >0$. We assume that 
\begin{equation*}
V_{t}(b)=g(t-b), \quad \quad V_{\mathbf{x}}(c)=f(\mathbf{x};c,\kappa )
\end{equation*}%
and that, for $\mathbf{u}=(u_{x},u_{y},u_{z})' \in $ $\mathds{S}^{2}$, $R_{(t,b)}=R_{\theta (t-b),\mathbf{u}}$, where $%
R_{\theta ,\mathbf{u}}$ is the rotation matrix of angle $\theta $
around an axis in the direction of $\mathbf{u}$. We have that
\begin{equation*}
R_{\theta ,\mathbf{u}}=\cos \theta I_{3}+\sin \theta \lbrack \mathbf{u}%
]_{\times }+(1-\cos \theta ) \Mb{u} \Mb{u}',
\end{equation*}%
where $I_{3}$ is the identity matrix of $\mathds{R}^{3}$ and $[\mathbf{u}]_{\times }$ the cross product matrix of $\mathbf{u}$, defined by
\begin{equation*}
\lbrack \mathbf{u}]_{\times }=
\begin{pmatrix}
0 & -u_{z} & u_{y} \\ 
u_{z} & 0 & -u_{x} \\ 
-u_{y} & u_{x} & 0%
\end{pmatrix}.
\end{equation*}

To the best of our knowledge, the resulting models are the first max-stable models on a sphere. Such models can of course be relevant in practice due to the natural spherical shape of planets and stars. Moreover, as before, this type of model has a max-autoregressive representation for an appropriate function $g$.

\medskip

\noindent \textbf{iii) Models of type 3: Schlather's representation with }$%
\mathcal{X}=\mathds{R}^{2}$

\medskip

For $d\in \mathds{N} \backslash \{ 0 \}$, let $\mathcal{C}_{d}=\mathcal{C}\left( 
\mathds{R}^{d},\mathds{R}_{+} \backslash \{ 0 \} \right) $ be the space of continuous functions
from $\mathds{R}^{d}$ to $\mathds{R}_{+} \backslash \{ 0 \}$. For this sub-class of models, $%
(E_{1},\mathcal{E}_{1},\mu _{1})$ and $(E_{2},\mathcal{E}_{2},\mu _{2})$ are probability spaces with $E_{1}=\mathcal{C}_{1}$, $E_{2}=\mathcal{C}_{2}$
and $\mu _{1}$ and $\mu _{2}$ are probability measures on $%
E_{1}$ and $E_{2}$, respectively. The function $V_{t}$ (respectively $V_{\mathbf{x}}$) is defined as the
natural projection from $\mathcal{C}_{1}$ (respectively $\mathcal{C}_{2}$) to $%
\mathds{R}_+$ such that 
\begin{equation*}
V_{t}(b)=b(t)\quad \text{and}\quad V_{\mathbf{x}}(c)=c(\mathbf{x}),
\end{equation*}%
with the conditions that $\mathds{E}[b_{i}(t)]=1$ and $\mathds{E}[c_{i}(%
\mathbf{x})]=1$. Note that for notational consistency, we use small letters for the stochastic processes $b$ and $c$. The spectral processes $c_{i}$ are assumed to be either
stationary and in this case $R_{(t,b)}\mathbf{x}=\mathbf{x}-t%
\boldsymbol{\tau }$ where $\boldsymbol{\tau \in }\mathds{R}^{2}$, or to be isotropic and in this case $R_{(t,b)}\mathbf{x}=A^{t}\mathbf{x}$ where $A$
is an orthogonal matrix ($R_{(t,b)}$ corresponds to a rotation).

\medskip

\noindent \textbf{iv) Models of type 4: Mixed representation with }$\mathcal{%
X}=\mathds{R}^{2}$

\medskip

We choose $E_{1}=\mathcal{I}$, $\mu
_{1}=\delta $, $\mathcal{X}=\mathds{R}^{2}$, $E_{2}=\mathcal{C}_{2}$. Let $g$
be a probability density function (case $\mathcal{I}=\mathds{R}$) or a discrete
probability distribution (case $\mathcal{I}=\mathds{Z}$) and $\mu _{2}$ a
probability measure on $\mathcal{C}_{2}$. We take
\begin{equation*}
V_{t}(b)=g(t-b) \quad \mbox{ and } \quad V_{\mathbf{x}}(c)=c(\mathbf{x})\text{.}
\end{equation*}%
As in the previous case, $V_{\mathbf{x}}$ is the natural projection from $\mathcal{C}_{2}$ to $\mathds{R}_{+}$. Once again, note that we use a small letter for the stochastic process $c$.
The spectral processes $c_{i}$ are assumed to be stationary and $%
R_{(t,b)}\mathbf{x}=\mathbf{x}-(t-b)\boldsymbol{\tau }$, where $\boldsymbol{\tau \in }%
\mathds{R}^{2}$.
As for models of types 1 and 2, the processes of this type can be written under a max-autoregressive form for a well-chosen function $g$.

\medskip

We now focus on the stationary distributions in space i.e. when we consider
a fixed time $t$ and look at the spatial dimension. For a fixed $t\in 
\mathcal{I}$, we define the process 
$\left( X_{t}(\mathbf{x})\right) _{\mathbf{x}\in \mathcal{X}}=\left( X(t,%
\mathbf{x})\right) _{\mathbf{x}\in \mathcal{X}}$. For two processes, $\overset{d}{=}$ denotes equality in distribution for any
finite dimensional vectors of the two processes.

\begin{Th}[Stationary distributions in space]
\label{Prop_fixed_t} For a fixed $t\in \mathcal{I}$, assume that for each $%
M \in \mathds{N} \backslash \{ 0 \}$, $b \in E_1$, $\mathbf{x}_{1},\dots ,%
\mathbf{x}_{M} \in \mathcal{X}$ and $z_{1},\dots
,z_{M}>0$, we have that
\begin{equation}
\int_{E_{2}}\bigvee_{m=1}^{M} \left \{ \frac{V_{R_{(t,b)}\mathbf{x}_{m}}(c)}{z_{m}} \right \} \mu
_{2}(dc)=\int_{E_{2}}\bigvee_{m=1}^{M} \left \{  \frac{V_{\mathbf{x}_{m}}(c)}{z_{m}} \right \} \mu
_{2}(dc).  \label{AssumptionStationarySpaceb}
\end{equation}%
Then we have that
\begin{equation*}
X_{t}(\mathbf{x})\overset{d}{=}\bigvee_{i=1}^{\infty } \left \{ U_{i}V_{\mathbf{x}%
}(C_{i}) \right \} ,\quad \mathbf{x}\in \mathds{R}^{2},
\end{equation*}%
where $(U_{i},C_{i})_{i\geq 1}$ are the points of a Poisson point process on $(0,\infty
)\times E_{2}$ with intensity $u^{-2}du\times \mu _{2}(dc)$.
Moreover, Assumption $\left( \ref{AssumptionStationarySpaceb}\right) $ is satisfied for models of types 1, 2, 3 and 4.
\end{Th}
We see in Theorem \ref{Prop_fixed_t} that the spectral separability and
the use of specific operators $R$ make the spectral function $V_{\mathbf{x}}$
(with its associated point process $(C_{i})_{i\geq 1}$) the function which
appears in the spatial spectral representation. The stationary distribution in space only depends on the spatial parameters of the model. This property is
interesting from a statistical point of view since any estimation procedure
can be simplified by considering in a first step the spatial parameters only
without taking into account the temporal ones (see Section 4). Note that the idea of using a transformation of space in \eqref{AssumptionStationarySpaceb} can also be found in \cite{strokorb2015tail}, in a different context.

\medskip

We now look at the marginal distributions in time i.e. when we consider a
fixed site $\mathbf{x}\in \mathcal{X}$. As previously, we define the process $\left( X_{\mathbf{x}}(t)\right) _{t\in \mathcal{I}}=\left( X(t,%
\mathbf{x})\right) _{t\in \mathcal{I}}.$

\begin{Th}[Marginal distributions in time]
\label{Prop_fixed_x} For a fixed $\mathbf{x}\in \mathcal{X}$, assume that
there exist two operators $S$ and $G$ from $\mathcal{X}$ to $\mathcal{X}$ such that 
\begin{equation}
R_{(t,b)}S_{(t)}\mathbf{x}=G_{(b)}\mathbf{x}.
\label{AssumptionStationaryTime}
\end{equation}%
Then we have that
\begin{equation*}
X_{S_{(t)}\mathbf{x}}(t)\overset{d}{=}\bigvee_{i=1}^{\infty
} \left \{ U_{i}V_{t}(B_{i}) \right \}, \quad t \in \mathds{Z},
\end{equation*}%
where $(U_{i},B_{i})_{i\geq 1}$ are the points of a Poisson point process on $%
(0,\infty )\times E_{1}$ with intensity $u^{-2}du\times \mu _{1}(db)$.
Assumption $\left( \ref{AssumptionStationaryTime}\right) $ is satisfied for
models of types 1 and 4 with $S_{(t)}\mathbf{x}=\mathbf{x}+t\boldsymbol{\tau }$, for models of type 2 with $S_{(t)}\mathbf{x}=R_{-\theta t,\mathbf{u}}\mathbf{x}$ and for models of type 3 with $S_{(t)}\mathbf{x}=\mathbf{x}+t\boldsymbol{\tau }$
or $S_{(t)}\mathbf{x}=A^{-t}\mathbf{x}$, where $\boldsymbol{\tau} \in \mathds{R}^2$ and $A$ is an orthogonal matrix.
\end{Th}
Contrary to Theorem \ref{Prop_fixed_t}, it is not possible to say that
the marginal distributions in time are those given by the temporal spectral
representation with the spectral function $V_{t}$ and its associated point
process $(B_{i})_{i\geq 1}$. In order to obtain such a representation, it is necessary to apply a time transformation $S_{(t)}$ on $\mathbf{x}$. As a consequence, it is difficult to estimate the temporal parameters separately since this
transformation is not necessarily known in practice. The transformation indeed depends on the type of model and the parameters we want to estimate. Note that if $R_{t,b}$ does not depend on $t$ (for instance the translation with $\boldsymbol{\tau}=\Mb{0}$), i.e. if space and time are fully separated in the spectral representation, then $S_{(t)}$ is equal to the identity.

\section{Markovian cases}

In this section, in the
case $\mathcal{I}=\mathds{R}$, $g$ is the density of a standard exponential
random variable whereas in the case $\mathcal{I}=\mathds{Z}$, $g$
corresponds to the probability weights of a geometric random variable: 
\begin{equation}
g(t)=\left\{ 
\begin{array}{ll}
\nu \exp (-\nu t)\ \mathds{I}_{\{t\geq 0\}} & \mbox{if }\mathcal{I}=\mathds{R%
} \\ 
(1-\phi )\phi ^{t}\ \mathds{I}_{\{t\geq 0\}} & \mbox{if }\mathcal{I}=\mathds{%
Z}%
\end{array}%
,\right.   \label{Eq_Function_g}
\end{equation}%
where $\nu >0$ and $\phi \in (0,1)$. 
We first consider models of type 1 and type 4 and then models of type 2.
The choice of the function $g$ in \eqref{Eq_Function_g} makes these spatio-temporal max-stable models Markovian.

\subsection{Markovian models of type 1 and type 4}
\label{Subsec_Markov_Types_1_4}

Recall that we assume the transformations $%
R_{(t,b)}$ to be translations: $R_{(t,b)}(\mathbf{x})=\mathbf{x}-(t-b)\boldsymbol{\tau }$, where $\boldsymbol{%
\tau }\in \mathds{R}^{2}$. The parameter $\boldsymbol{\tau }$ gives a
preferred direction of propagation of the process. In this context, we obtain
\begin{equation}
X(t,\mathbf{x})=\left\{ 
\begin{array}{ll}
\bigvee_{i\geq 1}\left\{ U_{i}\nu \exp (-\nu (t-B_{i}) ) \mathds{I}%
_{\{t-B_{i}\geq 0\}}V_{\mathbf{x}-(t-B_{i})\boldsymbol{\tau }%
}(C_{i})\right\}  & \mbox{if }\mathcal{I}=\mathds{R} \\ 
\bigvee_{i\geq 1}\left\{ U_{i}\phi (1-\phi )^{t-B_{i}}\mathds{I}%
_{\{t-B_{i}\geq 0\}}V_{\mathbf{x}-(t-B_{i})\boldsymbol{\tau }%
}(C_{i})\right\}  & \mbox{if }\mathcal{I}=\mathds{Z}%
\end{array}%
.\right.   \label{Model_deHaan_Special_Case}
\end{equation}%
Note that for $\mathcal{I}=\mathds{R}$, the function $g$ has been introduced
by \cite{dombry2014stationary}, under the form $g(t)=-\log (a)a^{t}\mathds{I}%
_{\{t\geq 0\}}$ for $a\in (0,1)$, in order to build the continuous-time
version of the real-valued max-AR(1) process.

The following result shows in particular that the process $X$ defined in %
\eqref{Model_deHaan_Special_Case} satisfies a stochastic recurrence
equation. Let us denote by $a$ the constant $\exp (-\nu )$ if $\mathcal{I}=%
\mathds{R}$ and the constant $\phi $ if $\mathcal{I}=\mathds{Z}$.

\begin{Th}
\label{Prop_Iteration} i) For all $t,s\in \mathcal{I}$ such that $s<t$, we have
that 
\begin{equation}
X(t,\mathbf{x})=\max (a^{s}X(t-s,\mathbf{x}-s\boldsymbol{\tau }%
),(1-a^{s})Z(t,\mathbf{x})),  \label{Eq_Iteration}
\end{equation}%
where the process $(Z(t,\mathbf{x}))_{\mathbf{x}\in \mathds{R}^{2}}$ is
independent of $(X(t-s,\mathbf{x}))_{\mathbf{x}\in \mathds{R}^{2}}$ and 
\begin{equation}
Z(t,\mathbf{x})\overset{d}{=}\bigvee_{i=1}^{\infty } \left \{ U_{i}V_{\mathbf{x}%
}(C_{i}) \right \} ,\quad \mathbf{x}\in \mathds{R}^{2},  
\label{SpectralrepresentationZ}
\end{equation}%
with $(U_{i},C_{i})_{i\geq 1}$ the points of a Poisson point process on $(0,\infty
)\times E_{2}$ of intensity $u^{-2}du\times \mu _{2}(dc)$. 

ii) Let $\mathcal{I}=\mathds{Z}$ and $((Z(t,\mathbf{x}))_{\mathbf{x}\in \mathds{R}^{2}})_{t\in \mathcal{I}}$
be a family of iid max-stable processes with spectral representation $\left( %
\ref{SpectralrepresentationZ}\right) $, we have that 
\begin{equation}
\label{Eq_Max_Integral_Representation}
X(t,\mathbf{x})\overset{d}{=}\bigvee_{j=0}^{\infty} \left\{ a^{j}(1-a)Z(t-j,%
\mathbf{x}-j\boldsymbol{\tau })\right\} .
\end{equation}
\end{Th}

From Theorem \ref{Prop_Iteration}, it can be seen that our model $X$
extends the real-valued MARMA$(1,0)$ process of \cite{davis1989basic} to the
spatial setting. The parameter $a$ measures the influence of the past,
whereas the parameter $\boldsymbol{\tau }$ represents some kind of specific
direction of propagation (contagion) in space. For the sake of ease of interpretation, consider the case where $%
\mathcal{I}=\mathds{Z}$ and $s=1$. The value at location $\mathbf{x}$ and
time $t$ is either related to the value at location $\mathbf{x}-\boldsymbol{%
\tau }$ at time $t-1$ or to the value of another process (the
innovation), $Z$, that characterizes a new event happening at location $%
\mathbf{x}$. If the value at location $\mathbf{x}-\boldsymbol{\tau }$ and $a$ are 
large, it is likely that there will be a propagation from location $\mathbf{x}-\boldsymbol{%
\tau }$ to location $\mathbf{x}$, i.e. contagion of the extremes, with an
attenuation effect. Contrary to the existing spatio-temporal max-stable
models, the dynamics is described by an equation that can be physically
interpreted. Note that the translation by the vector $-\boldsymbol{\tau }$
is one of the easiest transformations that allows to broaden the direct
extension of the real-valued MARMA$(1,0)$ model to a spatial setting.

Moreover, the combination of Theorems \ref{Prop_fixed_t} and \ref{Prop_Iteration}
shows that the stationary spatial distribution of the Markov
process/chain $\left( (X(t,\mathbf{x}) )_{\mathbf{x}\in \mathds{R}%
^{2}}\right) _{t\in \mathcal{I}}$ is the same as that of $Z$. It is important to remark that $\mathcal{C}_2$ (its state space) equipped e.g. with the topology induced by the distance $d(f_1,f_2)=\sup(\min(|f_1(\Mb{x})-f_2(\Mb{x})|,1): \Mb{x} \in \mathds{R}^2)$, for two functions $f_1$ and $f_2$, is not locally compact. Therefore, the theory developed e.g. in \cite%
{meyn2009markov} cannot be used to derive additional properties.

\medskip

We now consider the special case
\begin{equation}
X(t,\mathbf{x})=\max (aX(t-1,\mathbf{x}-\boldsymbol{\tau }),(1-a)Z(t,\mathbf{%
x})), \quad t \in \mathds{Z},
 \label{Eq_MARMA_01}
\end{equation}%
where $a\in (0,1)$, $\boldsymbol{\tau }\in \mathds{R}^{2}$ are both fixed
and $((Z(t,\mathbf{x}))_{\mathbf{x}\in \mathds{R}^{2}} )_{t\in \mathds{Z}}$
is a sequence of iid spatial max-stable processes with spectral
representation $\left( \ref{SpectralrepresentationZ}\right) $. The
general distribution function of this process is given in the next proposition.
\begin{Prop}
\label{Prop_General_Distribution_Function} For $M \in \mathds{N}\backslash
\{0\},t_{1}\leq \dots \leq t_{M}\in \mathds{Z}$, $\mathbf{x}_{1},\dots ,%
\mathbf{x}_{M} \in \mathds{R}^{2}$ and $z_{1},\dots ,z_{M} >0 $, we have that 
\begin{align}
& -\log \left( \mathds{P}(X(t_{1},\mathbf{x}_{1})\leq z_{1},\dots ,X(t_{M},%
\mathbf{x}_{M})\leq z_{M})\right)   \notag \\
& =\mathcal{V}_{\mathbf{x}_{1},\mathbf{x}_{2}-(t_{2}-t_{1})\boldsymbol{\tau 
},\dots ,\mathbf{x}_{M}-(t_{M}-t_{1})\boldsymbol{\tau }}\left( z_{1},\frac{%
z_{2}}{a^{t_{2}-t_{1}}},\dots ,\frac{z_{M}}{a^{t_{M}-t_{1}}}\right)   \notag
\\
& \ \ \ \ \ +\sum_{m=2}^{M-1}(1-a^{t_{m}-t_{m-1}})\mathcal{V}_{\mathbf{x}%
_{m},\mathbf{x}_{m+1}-(t_{m+1}-t_{m})\boldsymbol{\tau },\dots ,\mathbf{x}%
_{M}-(t_{M}-t_{m})\boldsymbol{\tau }}\left( z_{m},\frac{z_{m+1}}{%
a^{t_{m+1}-t_{m}}},\dots ,\frac{z_{M}}{a^{t_{M}-t_{m}}}\right)  \nonumber
\\& \ \ \ \ \  +\frac{%
1-a^{t_{M}-t_{M-1}}}{z_{M}},  \label{Eq_General_Distribution_Function}
\end{align}%
where $\mathcal{V}$ is the exponent function characterizing the spatial
distribution, defined by 
\begin{equation*}
\mathcal{V}_{\mathbf{x}_{1},\mathbf{x}_{2},\dots ,\mathbf{x}_{M}}\left(
z_{1},z_{2},\dots ,z_{M} \right) =\int_{\mathds{R}^2} \bigvee_{m=1}^{M} \left \{ \frac{V_{%
\mathbf{x}_{m}}(c)}{z_{m}} \right \} \mu_{2}(dc).
\end{equation*}
\end{Prop}

By using the approach developed by \cite{bienvenue2014likelihood}, the right-hand term of \eqref{Eq_General_Distribution_Function} can easily be
 computed provided that the distribution of $\left( V_{\mathbf{x}%
_{m}}(c)\right) _{m=1,\ldots ,M}$ $\ $with $c \sim \mu _{2}$ is absolutely
continuous with respect to the Lebesgue measure. This is the case for
example for the spatial Schlather and Brown-Resnick processes. 

It is easily shown that the models of types 1 and 4 are stationary in space and time.
In order to measure the spatio-temporal dependence, we propose extensions to the spatio-temporal
setting of  quantities that have been introduced in the spatial context. The first one is the spatio-temporal extremal coefficient function, stemming from the spatial version by \cite{schlather2003dependence}, which is
defined for all $t_{1},t_{2}\in \mathds{Z}$ and $%
\mathbf{x}_{1},\mathbf{x}_{2}\in \mathds{R}^{2}$ by 
\begin{equation*}
\mathds{P}(X(t_{1},\mathbf{x}_{1})\leq u,X(t_{2},\mathbf{x}_{2})\leq u)=\exp
\left( -\frac{\theta (t_{2}-t_{1},\mathbf{x}_{2}-\mathbf{x}_{1})}{u}\right) 
\text{,}\quad u>0.
\end{equation*}%
The second one is the spatio-temporal $\Phi _{1}$-madogram, coming from the spatial version introduced by \cite%
{cooley2006variograms}, where $\Phi _{1}$ is the standard Fr\'{e}chet probability distribution
function. It is defined by
\begin{equation*}
\nu _{\Phi _{1}}(t_{2}-t_{1},\mathbf{x}_{2}-\mathbf{x}_{1})=\frac{1}{2}%
\mathds{E}[|\Phi _{1}(X(t_{2},\mathbf{x}_{2}))-\Phi _{1}(X(t_{1},\mathbf{x}%
_{1}))|].
\end{equation*}%
\begin{Prop}
\label{Prop_Spatio_Temporal_Extremal_Coefficient} In the case of
\eqref{Eq_MARMA_01}, for $t_{1},t_{2}\in \mathds{Z}$ and $\mathbf{x}_{1},\mathbf{x}_{2}\in \mathds{R}%
^{2}$, the spatio-temporal extremal coefficient is given by 
\begin{equation}
\theta (t_{2}-t_{1},\mathbf{x}_{2}-\mathbf{x}_{1})=\mathcal{V}_{\mathbf{x}%
_{1},\mathbf{x}_{2}-(t_{2}-t_{1})\boldsymbol{\tau }}\left(
1,a^{t_{1}-t_{2}}\right) +1-a^{t_{2}-t_{1}}
\label{Eq_Spatio_Temporal_Extremal_Coefficient}
\end{equation}%
and the spatio-temporal $\Phi _{1}$-madogram of $X$ is given by 
\begin{equation}
\nu _{\Phi _{1}}(t_{2}-t_{1},\mathbf{x}_{2}-\mathbf{x}_{1})=\frac{1}{2}\frac{%
\theta (t_{2}-t_{1},\mathbf{x}_{2}-\mathbf{x}_{1})-1}{\theta (t_{2}-t_{1},%
\mathbf{x}_{2}-\mathbf{x}_{1})+1}=\frac{1}{2}-\frac{1}{\mathcal{V}_{\mathbf{x%
}_{1},\mathbf{x}_{2}-(t_{2}-t_{1})\boldsymbol{\tau }}\left(
1,a^{t_{1}-t_{2}}\right) +2-a^{t_{2}-t_{1}}}.  \label{Eq_F_Madogram}
\end{equation}%
\end{Prop}
Similarly, it would also be possible to extend the $\lambda $-madogram, introduced by \cite{naveau2009modelling}, to the spatio-temporal setting. 

Proposition \ref{Prop_Spatio_Temporal_Extremal_Coefficient} shows that we do
not fully separate space and time in the extremal dependence measure given by the extremal coefficient, even if $\boldsymbol{\tau }=0$. On the other hand, in the latter case, space and time are entirely separated in the spectral representation: $V_t$ only depends on time and $V_{R(t,b) \Mb{x}}=V_{\Mb{x}}$ only depends on space.

Furthermore, denoting $l=t_{2}-t_{1}$ and $\mathbf{h}=\mathbf{x}_{2}-\mathbf{x}_{1}$, we
have $\lim_{l \rightarrow \infty }\theta (l,\mathbf{h})=2$, showing
asymptotic time-independence. Moreover, from \cite{kabluchko2010ergodic}, we
deduce that,\ for a fixed $\mathbf{x}$, the process $\left( X_{\mathbf{x}%
}(t)\right) _{t\in \mathcal{I}}$ is strongly mixing in time. If $t_{1}=t_{2}$, $\lim_{\Vert \mathbf{h}\Vert
\rightarrow \infty }\theta (l,\mathbf{h})=2$ if and only if $X$ is strongly mixing
in space.

\medskip

Before showing some simulations, let us define the spatial Smith and Schlather models. Let $(U_i, C_i)_{i \geq 1}$ be the points of a Poisson point process on $(0, \infty) \times \mathds{R}^2$ with intensity $u^{-2} du \times \lambda_2(dc)$ and let $h_{\Sigma}$ denote the bivariate Gaussian density with mean $\Mb{0}$ and covariance matrix $\Sigma$. Then the spatial Smith model \citep{smith1990max} is defined as
$Z(\Mb{x}) =\bigvee_{i=1}^{\infty} \{ U_i h_{\Sigma}(\Mb{x}-C_i) \},$ for $\Mb{x} \in \mathds{R}^2$.
Let $(U_i)_{i \geq 1}$ be the points of a Poisson point process on $(0, \infty)$ with intensity $u^{-2} du$ and $Y_1, Y_2, \dots$ independent replications of the stochastic process $Y(\Mb{x})=\sqrt{2 \pi} \varepsilon(\Mb{x})$, for $\Mb{x} \in \mathds{R}^2$, where $\varepsilon$ is a stationary standard Gaussian process with correlation function $\rho(.)$. Then the spatial Schlather process \citep{schlather2002models} is defined as
$Z(\Mb{x})=\bigvee_{i=1}^{\infty} \{ U_i Y_i(\Mb{x}) \}$, for $\Mb{x} \in \mathds{R}^2$.

\begin{figure}[H]
\begin{minipage}[c]{-.3 \linewidth}
    \includegraphics[scale=0.45]{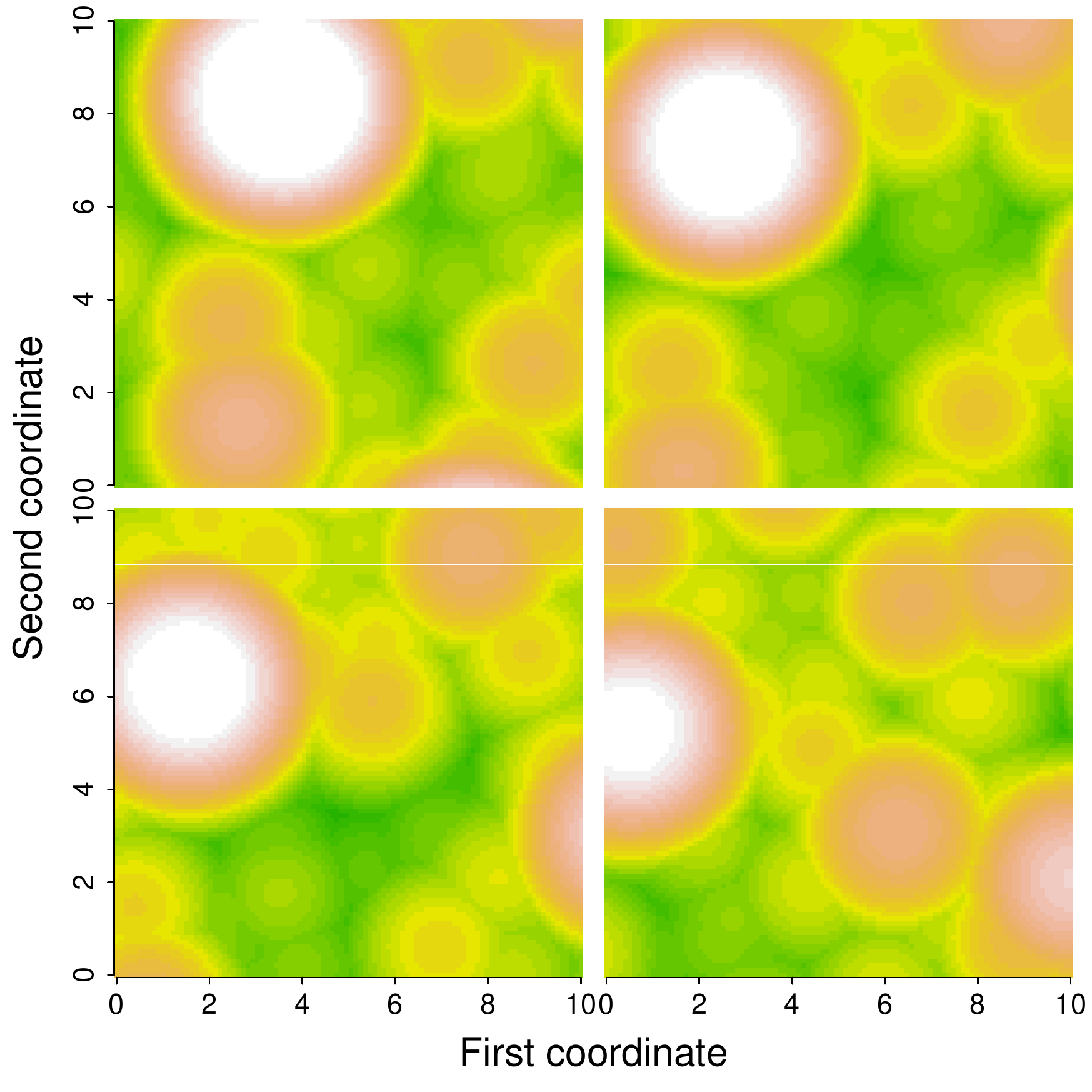} 
\end{minipage} \hfill
\begin{minipage}[c]{.47\linewidth}
    \includegraphics[scale=0.45]{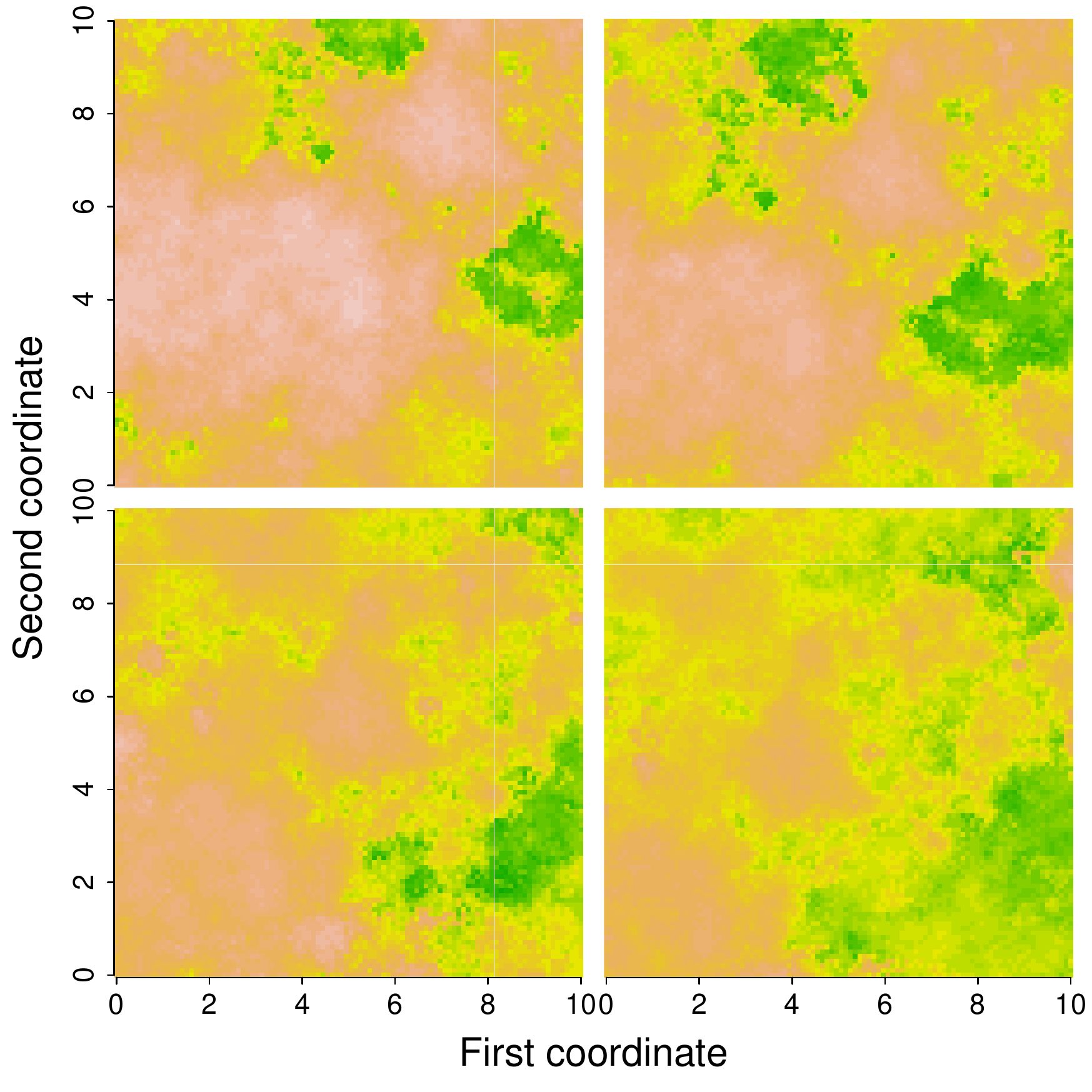} 
\end{minipage}
\caption{Simulation of the process \eqref{Eq_MARMA_01} for $Z$
being a spatial Smith process (left panel) and a spatial Schlather process (right panel). We depict the logarithm of the value of the process
in space. In both cases, the evolution over 4 periods is represented (top left: $t=1$, top
right: $t=2$, bottom left: $t=3$, bottom right: $t=4$).}    \label{Fig_Simulation_MARMA_Smith_Schlather}
\end{figure}

In the left panel of Figure \ref{Fig_Simulation_MARMA_Smith_Schlather}, we show the evolution of the
process \eqref{Eq_MARMA_01} when $Z$ is a spatial Smith process with covariance matrix 
\begin{equation*}
\Sigma=%
\begin{pmatrix}
1 & 0 \\ 
0 & 1%
\end{pmatrix}%
,
\end{equation*}%
$a=0.7$ and $\boldsymbol{\tau }=(-1,-1)^{\prime }$ (translation to the
bottom left). 
In the right panel of Figure \ref{Fig_Simulation_MARMA_Smith_Schlather}, we show the evolution of the process \eqref{Eq_MARMA_01} when $Z$ is a spatial Schlather process with correlation function of
type powered exponential, defined, for all $h\geq 0$ by $\rho(h)= \exp \left[ - \left( \frac{h}{c_1} \right)^{c_2} \right]$ for $c_1 >0$ and $0 < c_2 < 2$, where $c_1$ and $c_2$ are the range and the smoothing parameters, respectively. We take $c_{1}=3$, $c_{2}=1$ and, as previously, $a=0.7$ and $\boldsymbol{\tau }%
=(-1,-1)'$. Note that the process \eqref{Eq_MARMA_01} with $Z$ being the spatial Smith and the spatial Schlather process corresponds to models of types 1 and 4, respectively.
In both cases, we observe a translation of the main spatial structures (the
\textquotedblleft storms\textquotedblright\ in the case of the spatial Smith model)
to the bottom left, hence highlighting the usefulness of models like \eqref{Eq_MARMA_01} for phenomena that propagate in space.

\subsection{Markovian models of type 2}

Let $\mathcal{C}_{\mathds{S}^{2}}=\mathcal{C}(\mathds{S}^{2},\mathds{R}_{+} \backslash \{ 0 \})$ be the space of continuous functions from $\mathds{S}^{2}$ to $\mathds{R}_+ \backslash \{ 0 \}$. For 
$h\in \mathcal{C}_{\mathds{S}^{2}}$, let $\left\Vert h\right\Vert _{\infty
}=\sup_{\mathbf{x}\in \mathds{S}^{2}}h(\mathbf{x})$ and $d(h_{1},h_{2})=%
\left\Vert h_{1}-h_{2}\right\Vert _{\infty }$; then $C_{\mathds{S}^{2}}$ is a
Polish space \citep[see e.g.][Theorem 4.19]{kechris2012classical}.

By using the same arguments as in Theorem \ref{Prop_Iteration}, we
deduce that the models of type 2 satisfy the following stochastic recurrence
equation:
\begin{equation}
\label{Eq_Recurrence_Model_Type_2}
X(t,\mathbf{x})=\max (aX(t-1,R_{\theta ,\mathbf{u}}\mathbf{x}),(1-a)Z(t,%
\mathbf{x})),
\end{equation}%
where the process $(Z(t,\mathbf{x}))_{\mathbf{x}\in \mathds{S}^{2}}$ is
independent of $(X(t-1,\mathbf{x}))_{\mathbf{x}\in \mathds{S}^{2}}$ and is
such that%
\begin{equation}
Z(t,\mathbf{x})\overset{d}{=}\bigvee_{i=1}^{\infty } \left \{ U_{i}f(\mathbf{x};%
\boldsymbol{\mu }_{i},\kappa ) \right \} ,\quad \mathbf{x}\in \mathds{S}^{2},
\label{DistributionZspatial}
\end{equation}%
with $(U_{i},\boldsymbol{\mu }_{i})_{i\geq 1}$ the points of a Poisson point process on $%
(0,\infty )\times \mathds{S}^{2}$ with intensity $u^{-2}du\times d\lambda _{%
\mathds{S}^2}$.
Therefore, $\left( ( X(t,\mathbf{x}) )_{\mathbf{x}\in \mathds{S}^{2}}\right)
_{t\in \mathds{Z}}$ is a Markov chain with state space $\mathcal{C}_{\mathds{%
S}^{2}}$. Its transition kernel is denoted by $\mathcal{P}$. By using non standard results for Markov chains on Polish spaces (\cite{hairer2010p}), it is
possible to show that the transition kernel converges towards a unique invariant
measure at an exponential rate. Note that the state space considered in Section \ref{Subsec_Markov_Types_1_4}, $\mathcal{C}_2$, is not separable and hence is not a Polish space, yielding that the results of \cite{hairer2010p} cannot be used in that case.

Let $L$ be a function from $\mathcal{C}_{\mathds{S}^{2}}$ to $[0,\infty )$
and let us introduce a weighted supremum norm on the space of functions from 
$\mathcal{C}_{\mathds{S}^{2}}$ to $\mathds{R}$ in the following way. For $%
\varphi :\mathcal{C}_{\mathds{S}^{2}}\mapsto \mathds{R}$, we define
\begin{equation*}
\left\Vert \varphi \right\Vert _{L}=\sup_{h\in \mathcal{C}_{\mathds{S}^{2}}}%
\frac{|\varphi (h)|}{1+L(h)}.
\end{equation*}%
Finally, for a distribution $\pi $ on $\mathcal{C}_{\mathds{S}^{2}}$ and $%
\varphi :\mathcal{C}_{\mathds{S}^{2}}\mapsto \mathds{R}$, let us denote $\pi
(\varphi )=\int_{\mathcal{C}_{\mathds{S}^{2}}}\varphi (h)d\pi (h)$.
By Theorem 3.6 in \cite{hairer2010p}, we deduce the following result
about the geometric ergodicity of the Markov chain $\left( ( X(t,\mathbf{x}%
) )_{\mathbf{x}\in \mathds{S}^{2}}\right) _{t\in \mathds{Z}}$.

\begin{Th}[Geometric ergodicity]
\label{Prop_Geometric_Ergodicity} There exists a unique invariant measure $%
\pi $ for the Markov chain $\left( ( X(t,\mathbf{x}) )_{\mathbf{x}\in 
\mathds{S}^{2}}\right) _{t\in \mathds{Z}}$. Let, for $h\in \mathcal{C}_{%
\mathds{S}^{2}}$, $L(h)=\left\Vert h^{\gamma }\right\Vert _{\infty }$ with $%
0<\gamma <1$.\ There exist constants $C>0$ and $\rho \in (0,1)$ such that%
\begin{equation*}
\left\Vert \mathcal{P}^{n}\varphi -\pi (\varphi )\right\Vert _{L} \leq C\rho
^{n}\left\Vert \varphi -\pi (\varphi )\right\Vert _{L}
\end{equation*}%
holds for every measurable function $\varphi :\mathcal{C}_{\mathds{S}%
^{2}}\mapsto \mathds{R}$ such that $\left\Vert \varphi \right\Vert
_{L}<\infty $.
\end{Th}

\section{Estimation on simulated data}

In this section, we briefly discuss statistical inference for the process $
\eqref{Eq_MARMA_01}$. We denote by $\boldsymbol{\theta }$ the vector
gathering the parameters to be estimated. One possible method of estimation consists in using the pairwise likelihood \citep[see e.g.][]{davis2013statistical}, which requires the
knowledge of the bivariate density function for each $t_{1},t_{2}\in 
\mathds{R}$ and $\mathbf{x}_{1},\mathbf{x}_{2}\in \mathds{R}^{2}$. The latter is given in the following proposition.

\begin{Prop}
\label{Prop_Bivariate_Density_MARMA_01} For $t_{1},t_{2}\in \mathds{R}$, $%
\mathbf{x}_{1},\mathbf{x}_{2}\in \mathds{R}^{2}$ and $z_{1},z_{2}>0$, the bivariate density of the process \eqref{Eq_MARMA_01}
is given by 
\begin{align}
& f_{(t_{1},\mathbf{x}_{1}),(t_{2},\mathbf{x}_{2})}(z_{1},z_{2},\boldsymbol{%
\theta })  \notag \\
& =\exp \left( -\mathcal{V}_{\mathbf{x}_{1},\mathbf{x}_{2}-(t_{2}-t_{1})%
\boldsymbol{\tau }}\left( z_{1},\frac{z_{2}}{a^{t_{2}-t_{1}}}\right) -\frac{%
1-a^{t_{2}-t_{1}}}{z_{2}}\right)   \notag \\
& \ \ \ \times \Bigg[\left( -\frac{\partial }{\partial z_{1}}\mathcal{V}_{%
\mathbf{x}_{1},\mathbf{x}_{2}-(t_{2}-t_{1})\boldsymbol{\tau }}\left( z_{1},%
\frac{z_{2}}{a^{t_{2}-t_{1}}}\right) \right) \times \left( -\frac{\partial }{%
\partial z_{2}}\mathcal{V}_{\mathbf{x}_{1},\mathbf{x}_{2}-(t_{2}-t_{1})%
\boldsymbol{\tau }}\left( z_{1},\frac{z_{2}}{a^{t_{2}-t_{1}}}\right) +\frac{%
1-a^{t_{2}-t_{1}}}{z_{2}^{2}}\right)   \notag \\
& \ \ \ -\frac{\partial }{\partial z_{1}}\frac{\partial }{\partial z_{2}}%
\mathcal{V}_{\mathbf{x}_{1},\mathbf{x}_{2}-(t_{2}-t_{1})\boldsymbol{\tau }%
}\left( z_{1},\frac{z_{2}}{a^{t_{2}-t_{1}}}\right) \Bigg ].
\label{Eq_Bivariate_Density_MARMA_01}
\end{align}
\end{Prop}

Regarding the process $Z$ appearing in \eqref{Eq_MARMA_01}, we only consider the case of the spatial Smith model. Its covariance
matrix is denoted 
\begin{equation*}
\Sigma=%
\begin{pmatrix}
\sigma _{11} & \sigma _{12} \\ 
\sigma _{12} & \sigma _{22}%
\end{pmatrix}%
\end{equation*}%
and $\boldsymbol{\theta }$ is now given by $(\sigma _{11},\sigma
_{12},\sigma _{22},a,\boldsymbol{\tau }')'$. The bivariate
density function is given below.

\begin{Corr}
\label{Prop_Bivariate_Density_MARMA_01_Smith} We denote by $w_{1}=\frac{h_{1}%
}{2}+\frac{1}{h_{1}}\log \left( \frac{z_{2}}{a^{t_{2}-t_{1}}z_{1}}\right) $
and $v_{1}=h_{1}-w_{1}$, where 
\begin{equation*}
h_{1}=\sqrt{(\mathbf{x}_{2}-(t_{2}-t_{1})\boldsymbol{\tau }-\mathbf{x}%
_{1})^{\prime } \Sigma^{-1}(\mathbf{x}_{2}-(t_{2}-t_{1})\boldsymbol{%
\tau }-\mathbf{x}_{1})}.
\end{equation*}%
Let $t_{1},t_{2}\in \mathds{R}$, $\mathbf{x}_{1},\mathbf{x}_{2}\in \mathds{R}%
^{2}$ and $z_{1},z_{2} >0$. The bivariate density
function of the process \eqref{Eq_MARMA_01} with $Z$ being the spatial Smith process is given by 
\begin{align*}
& f_{(t_{1},\mathbf{x}_{1}),(t_{2},\mathbf{x}_{2})}(z_{1},z_{2},\boldsymbol{%
\theta }) \\
& =\exp \left( -\frac{\Phi (w_{1})}{z_{1}}-\frac{a^{t_{2}-t_{1}}\Phi (v_{1})%
}{z_{2}}-\frac{1-a^{t_{2}-t_{1}}}{z_{2}}\right)  \\
& \ \ \ \times \Bigg[\left( \frac{\Phi (w_{1})}{z_{1}^{2}}+\frac{\phi (w_{1})%
}{h_{1}z_{1}^{2}}-\frac{a^{t_{2}-t_{1}}\phi (v_{1})}{h_{1}z_{1}z_{2}}\right)
\times \left( \frac{a^{t_{2}-t_{1}}\Phi (v_{1})}{z_{2}^{2}}+\frac{%
a^{t_{2}-t_{1}}\phi (v_{1})}{h_{1}z_{2}^{2}}-\frac{\phi (w_{1})}{%
h_{1}z_{1}z_{2}}+\frac{1-a^{t_{2}-t_{1}}}{z_{2}^{2}}\right)  \\
& \ \ \ +\frac{v_{1}\phi (w_{1})}{h_{1}^{2}z_{1}^{2}z_{2}}+\frac{%
a^{t_{2}-t_{1}}w_{1}\phi (v_{1})}{h_{1}^{2}z_{1}z_{2}^{2}}\Bigg ],
\end{align*}%
where $\Phi $ and $\phi $ are the probability distribution
function and the probability density function of a standard Gaussian
random variable.
\end{Corr}
Assume that we observe the process at $M$ locations $\mathbf{x}_{1},\dots ,%
\mathbf{x}_{M}$ and $N$ dates $t_{1},\dots ,t_{N}$. Then, the spatio-temporal pairwise
log-likelihood is defined by \citep[see e.g.][]{davis2013statistical}
\begin{equation}
L_{P}^{ST}(\boldsymbol{\theta })=\sum_{i=1}^{N-1}\sum_{j=i+1}^{N}\sum_{k=1}^{M-1}%
\sum_{l=k+1}^{M}\omega _{i,j}\omega _{k,l}\log \left( f_{(t_{i},\mathbf{x}%
_{k}),(t_{j},\mathbf{x}_{l})}(z_{i,k},z_{j,l},\boldsymbol{\theta })\right) ,
\label{Eq_Pairwise_Likelihood}
\end{equation}%
where the $\omega _{i,j}$ and the $\omega _{k,l}$ are temporal and spatial
weights, respectively, and $z_{n,m}$ denotes the observation of the process at date $n$ and site $m$. Then, the maximum pairwise likelihood estimator is
given by$\boldsymbol{\ \hat{\theta}}=\mbox{argmax }L_{P}^{ST}(\boldsymbol{\theta }%
)$.

We will consider two different estimation schemes:
\begin{description}
\item[- ] \textbf{Scheme 1}:
As previously explained, due to Theorem \ref{Prop_fixed_t}, it is possible to separate the estimation of $\boldsymbol{\theta}_1=(\sigma _{11},\sigma
_{12},\sigma _{22})'$ and $\boldsymbol{\theta}_2=(a, \boldsymbol{\tau}')'$.
The estimation of $\boldsymbol{\theta}_1$ is carried out in a first step by maximizing the spatial pairwise log-likelihood \citep[see][]{padoan2010likelihood}.
Once $\boldsymbol{\theta}_1$ is known, we estimate $\boldsymbol{\theta}_2$ by maximizing $L_P^{ST}(\boldsymbol{\theta}_2)$ with respect to $\boldsymbol{\theta}_2$.
\item[- ] \textbf{Scheme 2}: We optimize $L_P^{ST}(\boldsymbol{\theta})$ with respect to $\boldsymbol{\theta}$, meaning that we estimate all parameters in a single step.
\end{description}

As illustration of the above, we simulate 100 times the process \eqref{Eq_MARMA_01} (where $Z$ is the spatial Smith process) with parameter $\boldsymbol{\theta }=(1,0,1,0.7,-1,-1)$ at $M$ sites and $N$ dates. We compute statistical summaries from the $100$ estimates obtained.
In both schemes, we optimize $L_P^{ST}$ with $\omega _{i,j}=1$ and $%
\omega _{k,l}=1$ for all $i=1,\dots ,N-1$, $j=i+1,\dots ,N$, $k=1,\dots ,M-1$
and $l=k+1,\dots ,M$. Tables \ref{Table_Scheme_1} and \ref{Table_Scheme_2} display the results for different
values of $M$ and $N$, in the cases of Scheme 1 and Scheme 2, respectively.
\begin{table}[h!]
\center
\resizebox{\textwidth}{!}{
\begin{tabular}{|c|c|c|c|c|c|c|c|}
\hline
\textbf{True} & \multicolumn{3}{|c|}{\textbf{Pairwise likelihood (M=20, N=20)}}& & \multicolumn{3}{|c|}{\textbf{Pairwise likelihood (M=30, N=30)}} \\ 
\hline\hline
\rowcolor{Gray} & Mean estimate & Mean bias & Stdev & & Mean estimate & Mean bias & Stdev \\ \hline
$\sigma_{11}$=1 & 1.139 & 0.139 & 0.421 & & 1.105 & 0.105 & 0.232 \\ \hline
$\sigma_{12}$=0 & 0.040 & 0.040 & 0.286 & & -0.024 & -0.024 & 0.162 \\ \hline
$\sigma_{22}$=1 & 1.185 & 0.185 & 0.325 & & 1.066 & 0.066 & 0.254 \\ \hline
a=0.7 & 0.707 & 0.007 & 0.059 & & 0.701 & 0.001 & 0.026 \\ \hline
$\boldsymbol{\tau}_1$=-1 & -0.990 & 0.010 & 0.123 & & -0.999 & 0.001 & 0.032 \\ \hline
$\boldsymbol{\tau}_2$=-1 & -0.990 & 0.010 & 0.101 & & -0.998 & 0.002 & 0.043  \\ \hline
\end{tabular}
}
\caption{Performance of the estimation in the case of Scheme 1. The mean estimate, the mean bias and the standard deviation are displayed.}
\label{Table_Scheme_1}
\end{table}

\begin{table}[h!]
\center
\resizebox{\textwidth}{!}{
\begin{tabular}{|c|c|c|c|c|c|c|c|}
\hline
\textbf{True} & \multicolumn{3}{|c|}{\textbf{Pairwise likelihood (M=20, N=20)}} & & \multicolumn{3}{|c|}{\textbf{Pairwise likelihood (M=30, N=30)}} \\ 
\hline\hline
\rowcolor{Gray} & Mean estimate & Mean bias & Stdev & & Mean estimate & Mean bias & Stdev \\ \hline
$\sigma_{11}$=1 & 1.288 & 0.288 & 0.678 & & 1.239 & 0.239 & 0.483 \\ \hline
$\sigma_{12}$=0 & 0.043 & 0.043 & 0.621 & & 0.057 & 0.057 & 0.314 \\ \hline
$\sigma_{22}$=1 & 1.453 & 0.453 & 1.159 & & 1.264 & 0.264 & 0.574 \\ \hline
a=0.7 & 0.706 & 0.006 & 0.050 & & 0.700 & 0.000 & 0.016 \\ \hline
$\boldsymbol{\tau}_1$=-1 & -0.998 & 0.002 & 0.115 & & -1.002 & -0.002 & 0.034 \\ \hline
$\boldsymbol{\tau}_2$=-1 & -0.982 & 0.018 & 0.111 & &  -1.003 & -0.003 & 0.035 \\ \hline
\end{tabular}
}
\caption{Performance of the estimation in the case of Scheme 2. The mean estimate, the mean bias and the standard deviation are displayed.}
\label{Table_Scheme_2}
\end{table}

For both schemes, the estimation is more accurate (the mean bias and the standard deviation decrease) as $M$ and $N$ increase.
Moreover, we observe that the estimation of the spatio-temporal parameters $a$ and $\boldsymbol{\tau}$ is satisfactory and clearly more accurate than that of the purely spatial parameters $\sigma_{11}, \sigma_{12}$ and $\sigma_{22}$ (the mean bias and the standard deviation are lower). Finally, the estimation of the purely spatial parameters is more accurate when using Scheme 1 (the mean bias and the standard deviation are lower). This stems probably from the fact that in Scheme 2, the number of pairs used is higher than in Scheme 1, introducing more variability. Indeed, contrary to what is assumed in the pairwise log-likelihood, the pairs considered are not independent. This dependence generates instability. For a discussion about the impact of the choice of pairs on estimation efficiency, see \cite{padoan2010likelihood}, p. 266 and 268. This finding shows that from a statistical point of view, spatio-temporal max-stable models that allow a separate estimation of the purely statistical parameters can be preferable; needless to say that a more extensive analysis would be needed at this point.

\section{Concluding remarks}

In order to overcome the defects of the spatio-temporal max-stable models
introduced in the literature, we propose a class of models where we partly
decouple the influence of time and space in the spectral representations.
Time has an influence on space through a bijective operator in space. Then, we propose several sub-classes of models where our operator is either a translation or a rotation. An advantage of the class of models we propose lies in the fact that it allows the roles of time and space to be distinct. Especially, the stationary distributions in space can differ from the marginal distributions in time. Moreover, the space operator allows to account for physical processes. Our models have both a continuous-time and a discrete-time version. 

Then, we consider a special case of some of our models where the function related to time in the spectral representation is the exponential density (continuous-time case) or takes as values the probabilities of a geometric random variable
(discrete-time case). In this context, the corresponding models become Markovian and have a useful max-autoregressive representation. They appear as an extension to a spatial setting of the real-valued MARMA$(1,0)$ process introduced by \cite{davis1989basic}. The main advantage of these models lies in the fact that the temporal
dynamics are explicit and easy to interpret. Moreover, these processes are
strongly mixing in time. We also show that the processes we introduce on the unit sphere of $\mathds{R}^3$ are geometrically ergodic. Finally, we briefly describe an 
inference method and show that it works well on simulated data, especially in the case of the parameters related to time. The detailed
study of possible estimation methodologies for our class of models will be considered in a subsequent paper.

\section*{Acknowledgements} 
Paul Embrechts acknowledges financial support by the Swiss Finance Institute (SFI). Erwan Koch would like to thank Andrea Gabrielli for interesting discussions. He also acknowledges RiskLab at ETH Zurich and the SFI for financial support. Christian Robert would like to
thank Mathieu Ribatet and Johan Segers for fruitful discussions on a related
topic when Johan Segers visited ISFA in October 2013.

\newpage
\appendix

\section{Poisson point process on $\mathds{Z}$}

\label{Sec_Poisson_Z}

For $A\subset \mathds{Z}$, let $\delta (A)=\sum_{z\in \mathds{Z}}\partial
_{\{z\}}(A)=\# A$ where $\# $ stand for the cardinality of a set. The point
process $N\left( A\right) =\sum_{k\in A}N_{k}$, $A\subset \mathds{Z}$,
defines an homogeneous Poisson point process on $\mathds{Z}$ with constant
intensity density function equal to $1$ since:

\begin{itemize}
\item due to the additivity of the Poisson distribution, $N\left(A\right)$
is Poisson distributed with parameter $\delta (A)$;

\item for any $l\geq 1$ and $A_{1},\ldots ,A_{l}$ disjoint sets in $\mathds{Z%
}$, the $N\left( A_{i}\right)$, $i=1, \dots, l,$ are independent random
variables.
\end{itemize}

\section{Proofs}

\subsection{For Theorem \protect\ref{Prop_fixed_t}}

\begin{proof}
For $M \in \mathds{N}\backslash \{0\}$, let $t\in \mathds{Z}$, $\mathbf{x}%
_{1},\dots ,\mathbf{x}_{M}\in \mathds{R}^{2}$ and $z_{1},\dots ,z_{M}>0$, we
deduce by \eqref{Eq_Finite_Dimensional_Distribution} and $\left( \ref{AssumptionStationarySpaceb}\right) $ that%
\begin{align*}
&-\log \left( \mathds{P}(X(t,\mathbf{x}_{1})\leq z_{1},\dots ,X(t,\mathbf{x}%
_{M})\leq z_{M})\right) 
\\& =\int_{E_{1}\times E_{2}}\bigvee_{m=1}^{M} \left \{ \frac{%
V_{t}(b)V_{R_{(t,b)}\mathbf{x}_{m}}(c)}{z_{m}} \right \} \mu _{1}(db)\mu _{2}(dc) 
\\&=\int_{E_{1}}V_{t}(b)\left( \int_{E_{2}}\bigvee_{m=1}^{M} \left \{ \frac{V_{R_{(t,b)}%
\mathbf{x}_{m}}(c)}{z_{m}} \right \} \mu _{2}(dc)\right) \mu _{1}(db) 
\\&=\int_{E_{2}}\bigvee_{m=1}^{M} \left \{ \frac{V_{\mathbf{x}_{m}}(c)}{z_{m}} \right \} \mu
_{2}(dc)\int_{E_{1}}V_{t}(b)\mu _{1}(db) 
\\&=\int_{E_{2}}\bigvee_{m=1}^{M} \left \{ \frac{V_{\mathbf{x}_{m}}(c)}{z_{m}} \right \} \mu
_{2}(dc).
\end{align*}
We now show that Assumption \eqref{AssumptionStationarySpaceb} is satisfied for models of types 1, 2, 3 and 4. \\
For models of type 1, we have $E_2=\mathds{R}^2$, $V_{R_{(t,b)}\mathbf{x}_{m}}(c)=f(R_{(t,b)}%
\mathbf{x}_{m}-c)=f(\Mb{x}_m-(c+(t-b) \boldsymbol{\tau}))=V_{\Mb{x}_m}(c+(t-b) \boldsymbol{\tau})$ and $\mu _{2}=\lambda _{2}$. Since $\lambda_2$ is invariant under
translation, we
derive by a change of variable that%
\begin{align*}
\int_{E_{2}}\bigvee_{m=1}^{M} \left \{ \frac{V_{R_{(t,b)}\mathbf{x}_{m}}(c)}{z_{m}} \right \} \mu
_{2}(dc) &=\int_{E_{2}}\bigvee_{m=1}^{M} \left \{ \frac{V_{\mathbf{x}_{m}}(c\mathbf{+}%
(t-b)\boldsymbol{\tau })}{z_{m}} \right \} \mu _{2}(dc)
\\& =\int_{E_{2}}\bigvee_{m=1}^{M}\left \{ \frac{%
V_{\mathbf{x}_{m}}(c)}{z_{m}} \right \} \mu _{2}(dc).
\end{align*}
For models of type 2, we have $E_2=\mathds{S}^2$ and
\begin{equation*}
V_{R_{(t,b)}\mathbf{x}_{m}}(c)=f(R_{\theta (t-b), \Mb{u}}\mathbf{x}_{m};c,\kappa )=%
\frac{\kappa }{4\pi \sinh \kappa }\exp \left( \kappa (R_{-\theta (t-b), \Mb{u}}c
\mathbf{)}' \mathbf{x}\right) = V_{\mathbf{x}_{m}}(R_{-\theta
(t-b), \Mb{u}}c)
\end{equation*}%
and it follows, since $\mu _{2}=\lambda _{\mathds{S}^2}$ is invariant under rotation, that%
\begin{align*}
\int_{E_{2}}\bigvee_{m=1}^{M} \left \{ \frac{V_{R_{(t,b)}\mathbf{x}_{m}}(c)}{z_{m}} \right \} \mu
_{2}(dc) &=\int_{E_{2}}\bigvee_{m=1}^{M} \left \{ \frac{V_{\mathbf{x}_{m}}(R_{-\theta
(t-b), \Mb{u}}c)}{z_{m}} \right \} \mu _{2}(dc)
\\& =\int_{E_{2}}\bigvee_{m=1}^{M}\left \{ \frac{V_{\mathbf{x}%
_{m}}(c)}{z_{m}} \right \} \mu _{2}(dc).
\end{align*}%
 \\
For models of type 3, we have $V_{\Mb{x}}(c)=c(\mathbf{x})$. Thus, if $R_{(t,b)}\mathbf{x}=\mathbf{x}-t\boldsymbol{\tau }$, we have $V_{R_{(t,b)}\mathbf{x}_{m}}(c)=c(
\mathbf{x}_m-t \boldsymbol{\tau})$. Thus, we  deduce
by stationarity that%
\begin{equation*}
\int_{E_{2}}\bigvee_{m=1}^{M} \left \{ \frac{V_{R_{(t,b)}\mathbf{x}_{m}}(c)}{z_{m}} \right \} \mu
_{2}(dc)
=\int_{E_{2}}\bigvee_{m=1}^{M} \left \{ \frac{V_{\mathbf{x}%
_{m}}(c)}{z_{m}} \right \} \mu _{2}(dc).
\end{equation*}%
If $R_{(t,b)}\mathbf{x}=A^{t}\mathbf{x}$, we have $V_{R_{(t,b)}\mathbf{x}_{m}}(c)=c(A^t \Mb{x}_m)$. Hence, we obtain by isotropy that%
\begin{equation*}
\int_{E_{2}}\bigvee_{m=1}^{M} \left \{ \frac{V_{R_{(t,b)}\mathbf{x}_{m}}(c)}{z_{m}} \right \} \mu
_{2}(dc)
=\int_{E_{2}}\bigvee_{m=1}^{M} \left \{ \frac{V_{\mathbf{x}_{m}}(c)}{z_{m}%
} \right \} \mu _{2}(dc).
\end{equation*}%
For models of type 4, we have $V_{\Mb{x}}(c)=c(\mathbf{x})$. Thus, if $R_{(t,b)}\mathbf{x}=\mathbf{x}-(t-b)\boldsymbol{\tau }$, we have $V_{R_{(t,b)}\mathbf{x}_{m}}(c)=c(\Mb{x}_m-(t-b) \boldsymbol{\tau})$. Hence, we
deduce by stationarity that%
\begin{equation*}
\int_{E_{2}}\bigvee_{m=1}^{M} \left \{ \frac{V_{R_{(t,b)}\mathbf{x}_{m}}(c)}{z_{m}} \right \} \mu
_{2}(dc)=\int_{E_{2}}\bigvee_{m=1}^{M} \left \{ \frac{V_{\mathbf{%
x}_{m}}(c)}{z_{m}} \right \} \mu _{2}(dc).
\end{equation*}
\end{proof}

\subsection{For Theorem \protect\ref{Prop_fixed_x}}

\begin{proof}
For $M \in \mathds{N}\backslash \{0\}$, $t_{1}, \dots, t_{M} \in 
\mathds{Z}$, $\mathbf{x}\in \mathds{R}^{2}$ and $z_{1},\dots ,z_{M}>0$, we
have that 
\begin{align*}
& -\log \left( \mathds{P}(X(t_{1},S_{(t_{1})}\mathbf{x})\leq z_{1},\dots
,X(t_{M},S_{(t_{M})}\mathbf{x})\leq z_{M})\right)  
\\&=\int_{E_{1}\times
E_{2}}\bigvee_{m=1}^{M} \left \{ \frac{V_{t_{m}}(b)V_{R_{(t_{m},b)}S_{(t_{m})}\mathbf{x%
}}(c)}{z_{m}} \right \} \mu _{1}(db)\mu _{2}(dc) 
\\&=\int_{E_{1}\times E_{2}}\bigvee_{m=1}^{M} \left \{ \frac{V_{t_{m}}(b)V_{G_{(b)}%
\mathbf{x}}(c)}{z_{m}} \right \} \mu _{1}(db)\mu _{2}(dc) 
\\&=\int_{E_{1}}\left( \bigvee_{m=1}^{M} \left \{ \frac{V_{t_{m}}(b)}{z_{m}} \right \}
\int_{E_{2}}V_{G_{(b)}\mathbf{x}}(c)\mu _{2}(dc)\right) \mu _{1}(db) 
\\&= \int_{E_{1}}\bigvee_{m=1}^{M} \left \{ \frac{V_{t_{m}}(b)}{z_{m}} \right \} \mu _{1}(db).
\end{align*}%
Moreover, it is easy to show that Assumption \eqref{AssumptionStationaryTime} is satisfied for models of types 1, 2, 3 and 4 with the operators $S_{(t)}$ that are given.
\end{proof}

\subsection{For Theorem \protect\ref{Prop_Iteration}}

\begin{proof}
i) Let us consider the case $\mathcal{I}=\mathds{R}$ (the case $\mathcal{%
I}=\mathds{Z}$ is similar). We have that
\begin{eqnarray*}
X(t,\mathbf{x}) &=&\bigvee_{i=1}^{\infty} \left\{ U_{i}\nu \exp (-\nu (t-B_{i}))%
\mathds{I}_{\{t-B_{i}\geq 0\}}V_{\mathbf{x}-(t-B_{i})\boldsymbol{\tau }%
}(C_{i})\right\}  \\
&=&\bigvee_{i=1}^{\infty} \left\{ U_{i}\nu \exp (-\nu (s+t-s-B_{i})) \mathds{I}%
_{\{s+t-s-B_{i}\geq 0\}}V_{\mathbf{x}-(s+t-s-B_{i})\boldsymbol{\tau }%
}(C_{i})\right\}  \\
&=&\max \Big(\exp (-\nu )^{s}\bigvee_{i=1}^{\infty} \left\{ U_{i} \nu \exp (-\nu
(t-s-B_{i}))\mathds{I}_{\{t-s-B_{i}\geq 0\}}V_{\mathbf{x}-s\boldsymbol{\tau }%
-(t-s-B_{i})\boldsymbol{\tau }}(C_{i}))\right\},  \\
&&\bigvee_{i=1}^{\infty} \left\{ U_{i} \nu \exp (-\nu (t-B_{i}))\mathds{I}%
_{\{t\geq B_{i}>t-s\}}V_{\mathbf{x}-(t-B_{i})\boldsymbol{\tau }%
}(C_{i}))\right\} \Big) \\
&=&\max \Big(\exp (-\nu )^{s}X(t-s,\mathbf{x}-s\boldsymbol{\tau }),(1-\exp
(-\nu )^{s})Z(t,\mathbf{x})\Big),
\end{eqnarray*}%
where%
\begin{equation*}
Z(t,\mathbf{x})=\frac{1}{(1-\exp (-\nu )^{s})}\bigvee_{i\geq 1}\left\{
U_{i} \nu \exp (-\nu (t-B_{i}))\mathds{I}_{\{t\geq B_{i}>t-s\}}V_{\mathbf{x}%
-(t-B_{i})\boldsymbol{\tau }}(C_{i}))\right\} .
\end{equation*}
Since the sets $\{t\geq B>t-s\}$ and $\{t-s\geq B\}$ are disjoint, the
Poisson point processes $\{(U_{i},B_{i},C_{i}),i:t\geq B_{i}>t-s\}$ and $%
\{(U_{i},B_{i},C_{i}),i:t-s\geq B_{i}\}$ are independent and it follows that 
$(X(t-s,\mathbf{x}))_{\Mb{x} \in \mathds{R}^2}$ and $(Z(t,\mathbf{x}))_{\Mb{x} \in \mathds{R}^2}$ are also
independent.

\medskip

\noindent We now show that 
\begin{equation*}
Z(t,\mathbf{x})\overset{d}{=}\bigvee_{i=1}^{\infty } \left \{ U_{i}V_{\mathbf{x}%
}(C_{i}) \right \},\quad \mathbf{x}\in \mathds{R}^{2},
\end{equation*}
where $(U_{i},C_{i})_{i\geq 1}$ are the points of a Poisson point process on $(0,\infty
)\times E_{2}$ of intensity $u^{-2}du\times \mu _{2}(dc)$.
Let $\left( U_{i},B_{i},C_{i}\right) _{i\geq 1}$ be the points of a Poisson point process
on $(0,\infty )\times \mathds{R}\times E_{2}$ with intensity $u^{-2}du\times
db\times \mu _{2}(dc)$. For $M\in \mathds{N}\backslash \{0\}$, let $\mathbf{x%
}_{1},\dots ,\mathbf{x}_{M}\in \mathds{R}^{2}$ and $z_{1},\dots ,z_{M}>0$.
We consider the set 
\begin{align*}
& B_{z_{1},\dots ,z_{M}}
\\& =\{(u,b,c):u\nu \exp (-\nu (t-b)) \mathds{I}_{\{t\geq
b>t-s\}}V_{\mathbf{x}_{m}-(t-b)\boldsymbol{\tau }}(c)>z_{m}%
\mbox{
for at least one }m=1,\dots ,M\}.
\end{align*}%
Denoting by $\bigwedge$ the min-operator, the Poisson measure of $B_{z_{1},\dots ,z_{M}}$ is%
\begin{eqnarray*}
\Lambda (B_{z_{1},\dots ,z_{M}}) &=&\int_{E_{2}}\int_{\mathds{R}%
}\int_{0}^{\infty }\mathds{I}_{\left\{ u> \bigwedge_{m=1}^M \left \{ \frac{z_{m}}{\nu \exp
(-\nu (t-b)) \mathds{I}_{\{t\geq b>t-s\}}V_{\mathbf{x}_{m}-(t-b)\boldsymbol{%
\tau }}(c)} \right \}  \right\} }u^{-2}du\times db \times \mu _{2}(dc) \\
&=&\int_{E_{2}}\int_{\mathds{R}}\bigvee_{m=1}^{M} \left \{ \frac{\nu \exp (-\nu (t-b))%
\mathds{I}_{\{t\geq b>t-s\}}V_{\mathbf{x}_{m}-(t-b)\boldsymbol{\tau }}(c)}{%
z_{m}} \right \} db\times \mu _{2}(dc) \\
&=&\nu \exp (-\nu t) \int_{\mathds{R}}\mathds{I}_{\{t\geq b>t-s\}}\exp (\nu
b) \left( \int_{E_{2}}\bigvee_{m=1}^{M} \left \{ \frac{V_{\mathbf{x}_{m}-(t-b)%
\boldsymbol{\tau }}(c)}{z_{m}} \right \} \mu _{2}(dc)\right) db.
\end{eqnarray*}%
Since 
\begin{equation*}
\int_{E_{2}}\bigvee_{m=1}^{M} \left \{ \frac{V_{\mathbf{x}_{m}-(t-b)\boldsymbol{\tau }%
}(c)}{z_{m}} \right \} \mu _{2}(dc)=\int_{E_{2}}\bigvee_{m=1}^{M}\left\{ \frac{V_{\mathbf{x}%
_{m}}(c)}{z_{m}} \right \} \mu _{2}(dc),
\end{equation*}%
we deduce that%
\begin{eqnarray*}
\Lambda (B_{z_{1},\dots ,z_{M}}) &=&\nu \exp (-\nu t) \left(
\int_{E_{2}}\bigvee_{m=1}^{M} \left \{ \frac{V_{\mathbf{x}_{m}}(c)}{z_{m}} \right \} \mu
_{2}(dc)\right) \int_{\mathds{R}}\mathds{I}_{\{t\geq b>t-s\}}\exp (\nu b) db
\\
&=&(1-\exp (-\nu )^{s})\left( \int_{E_{2}}\bigvee_{m=1}^{M} \left \{ \frac{V_{\mathbf{x%
}_{m}}(c)}{z_{m}} \right \} \mu _{2}(dc)\right) .
\end{eqnarray*}%
It follows that%
\begin{align*}
& -\log \left( \mathds{P}(Z(t,\mathbf{x}_{1})\leq z_{1},\dots ,Z(t,\mathbf{x}%
_{M})\leq z_{M})\right)
\\& = \Lambda (B_{(1-\exp (-\nu )^{s})z_{1},\dots
,(1-\exp (-\nu )^{s})z_{M}})=\int_{E_{2}}\bigvee_{m=1}^{M} \left \{ \frac{V_{\mathbf{x}%
_{m}}(c)}{z_{m}} \right \} \mu _{2}(dc).
\end{align*}

ii) It is easily shown that the right-hand side of \eqref{Eq_Max_Integral_Representation} is solution of \eqref{Eq_Iteration}. Moreover, as in \cite{davis1989basic}, this solution is unique, yielding \eqref{Eq_Max_Integral_Representation}.
\end{proof}

\subsection{For Proposition \protect\ref{Prop_General_Distribution_Function}}

\begin{proof}
For the sake of notational simplicity, we only give the proof in the case $%
M=3$; this proof can easily be extended. Using the independence of the
replications $(Z(t,\mathbf{x}))_{\mathbf{x}\in \mathds{R}^{2}}$ and changes
of indices, we obtain 
\begin{align}
& \mathds{P}\Bigg(\bigvee_{j=0}^{J}\left\{ a^{j}(1-a)Z(t_{1}-j,\mathbf{x}_{1}-j%
\boldsymbol{\tau })\right\} \leq z_{1},\bigvee_{j=0}^{J}\left\{
a^{j}(1-a)Z(t_{2}-j,\mathbf{x}_{2}-j\boldsymbol{\tau })\right\} \leq z_{2}, \nonumber
\\& \ \ \ \ \ \bigvee_{j=0}^{J}\left\{ a^{j}(1-a)Z(t_{3}-j,\mathbf{x}%
_{3}-j\boldsymbol{\tau })\right\} \leq z_{3}\Bigg)  \notag \\
& =\mathds{P}\Bigg(\bigvee_{j=0}^{J}\left\{ a^{j}(1-a)Z(t_{1}-j,\mathbf{x}_{1}-j%
\boldsymbol{\tau })\right\} \leq
z_{1}, \nonumber
\\& \ \ \ \ \ \ \bigvee_{j=t_{1}-t_{2}}^{J+t_{1}-t_{2}}\left\{
a^{j+t_{2}-t_{1}}(1-a)Z(t_{1}-j,\mathbf{x}_{2}-(j+t_{2}-t_{1})\boldsymbol{%
\tau }\right\} \leq z_{2},  \notag \\
& \ \ \ \ \ \ \bigvee_{j=t_{1}-t_{3}}^{J+t_{1}-t_{3}}\left\{
a^{j+t_{3}-t_{1}}(1-a)Z(t_{1}-j,\mathbf{x}_{3}-(j+t_{3}-t_{1})\boldsymbol{%
\tau }\right\} \leq z_{3}\Bigg)  \notag \\
& =\mathds{P}\Bigg(\bigvee_{j=0}^{J+t_{1}-t_{3}}\left\{ a^{j}(1-a)Z(t_{1}-j,%
\mathbf{x}_{1}-j\boldsymbol{\tau })\right\} \leq
z_{1}, \nonumber
\\& \ \ \ \ \ \ \ \ \ \bigvee_{j=0}^{J+t_{1}-t_{3}}\left\{ a^{j+t_{2}-t_{1}}(1-a)Z(t_{1}-j,%
\mathbf{x}_{2}-(j+t_{2}-t_{1})\boldsymbol{\tau }\right\} \leq z_{2},  \notag
\\
& \ \ \ \ \ \ \ \ \ \bigvee_{j=0}^{J+t_{1}-t_{3}}\left\{
a^{j+t_{3}-t_{1}}(1-a)Z(t_{1}-j,\mathbf{x}_{3}-(j+t_{3}-t_{1})\boldsymbol{%
\tau }\right\} \leq z_{3}\Bigg)  \notag \\
& \ \ \ \times \mathds{P}\Bigg(\bigvee_{j=t_{1}-t_{2}}^{-1}\left\{
a^{j+t_{2}-t_{1}}(1-a)Z(t_{1}-j,\mathbf{x}_{2}-(j+t_{2}-t_{1})\boldsymbol{%
\tau })\right\} \leq z_{2},  \notag \\
& \ \ \ \ \ \ \ \ \ \ \ \bigvee_{j=t_{1}-t_{2}}^{-1}\left\{
a^{j+t_{3}-t_{1}}(1-a)Z(t_{1}-j,\mathbf{x}_{3}-(j+t_{3}-t_{1})\boldsymbol{%
\tau })\right\} \leq z_{3}\Bigg)  \notag \\
& \ \ \ \times \mathds{P}\left( \bigvee_{j=t_{1}-t_{3}}^{t_{1}-t_{2}-1}\left\{
a^{j+t_{3}-t_{1}}(1-a)Z(t_{1}-j,\mathbf{x}_{3}-(j+t_{3}-t_{1})\boldsymbol{%
\tau })\right\} \leq z_{3}\right)   \notag \\
& \ \ \ \times \mathds{P}\left( \bigvee_{J+t_{1}-t_{2}+1}^{J}\left\{
a^{j}(1-a)Z(t_{1}-j,\mathbf{x}_{1}-j\boldsymbol{\tau })\right\} \leq
z_{1}\right)   \notag \\
& \ \ \ \times \mathds{P}\Bigg(\bigvee_{j=J+t_{1}-t_{3}+1}^{J+t_{1}-t_{2}}\left%
\{ a^{j}(1-a)Z(t_{1}-j,\mathbf{x}_{1}-j\boldsymbol{\tau })\right\} \leq
z_{1}),  \notag \\
& \ \ \ \ \ \ \ \ \ \ \  \bigvee_{j=J+t_{1}-t_{3}+1}^{J+t_{1}-t_{2}}\left\{
a^{j+t_{2}-t_{1}}(1-a)Z(t_{1}-j,\mathbf{x}_{2}-(j+t_{2}-t_{1})\boldsymbol{%
\tau })\right\} \leq z_{2}\Bigg).  \label{Eq_Distribution_Function_1}
\end{align}%
Using the independence of the replications $(Z(t,\mathbf{x}))_{\mathbf{x}\in 
\mathds{R}^{2}}$, the stationarity of the processes $(Z(t,\mathbf{x}%
))_{\mathbf{x}\in \mathds{R}^{2}}$ and the homogeneity of order $-1$ of $\mathcal{V}$,
we obtain 
\begin{align}
& \mathds{P}\Bigg( \bigvee_{j=0}^{J+t_{1}-t_{3}}\left\{ a^{j}(1-a)Z(t_{1}-j,%
\mathbf{x}_{1}-j\boldsymbol{\tau })\right\} \leq
z_{1}, \nonumber
\\& \ \ \ \ \  \bigvee_{j=0}^{J+t_{1}-t_{3}}\left\{ a^{j+t_{2}-t_{1}}(1-a)Z(t_{1}-j,%
\mathbf{x}_{2}-(j+t_{2}-t_{1})\boldsymbol{\tau })\right\} \leq z_{2},  \notag
\\
& \ \ \ \ \ \bigvee_{j=0}^{J+t_{1}-t_{3}}\left\{
a^{j+t_{3}-t_{1}}(1-a)Z(t_{1}-j,\mathbf{x}_{3}-(j+t_{3}-t_{1})\boldsymbol{%
\tau })\right\} \leq z_{3}\Bigg)  \notag \\
& =\prod_{j=0}^{J+t_{1}-t_{3}}\mathds{P}\Bigg(Z(t_{1}-j,\mathbf{x}_{1}-j%
\boldsymbol{\tau })\leq \frac{z_{1}}{a^{j}(1-a)},Z(t_{1}-j,\mathbf{x}%
_{2}-(j+t_{2}-t_{1})\boldsymbol{\tau })\leq \frac{z_{2}}{%
a^{j+t_{2}-t_{1}}(1-a)},  \notag \\
& \ \ \ \ \ \ Z(t_{1}-j,\mathbf{x}_{3}-(j+t_{3}-t_{1})\boldsymbol{\tau }\leq 
\frac{z_{3}}{a^{j+t_{3}-t_{1}}(1-a)}\Bigg)  \notag \\
& =\prod_{j=0}^{J+t_{1}-t_{3}}\exp \left( -\mathcal{V}_{\mathbf{x}_{1}-j%
\boldsymbol{\tau },\mathbf{x}_{2}-(t_{2}-t_{1})\boldsymbol{\tau }-j%
\boldsymbol{\tau },\mathbf{x}_{3}-(t_{3}-t_{1})\boldsymbol{\tau }-j%
\boldsymbol{\tau }}\left( \frac{z_{1}}{a^{j}(1-a)},\frac{z_{2}}{%
a^{j}a^{t_{2}-t_{1}}(1-a)},\frac{z_{3}}{a^{j}a^{t_{3}-t_{1}}(1-a)}\right)
\right)   \notag \\
& =\exp \left( - \mathcal{V}_{\mathbf{x}_{1},\mathbf{x}_{2}-(t_{2}-t_{1})%
\boldsymbol{\tau },\mathbf{x}_{3}-(t_{3}-t_{1})\boldsymbol{\tau }}\left(
z_{1},\frac{z_{2}}{a^{t_{2}-t_{1}}},\frac{z_{3}}{a^{t_{3}-t_{1}}}\right)
\times (1-a)\sum_{j=0}^{J+t_{1}-t_{3}}a^{j}\right)   \notag \\
& =\exp \left( - \mathcal{V}_{\mathbf{x}_{1},\mathbf{x}_{2}-(t_{2}-t_{1})%
\boldsymbol{\tau },\mathbf{x}_{3}-(t_{3}-t_{1})\boldsymbol{\tau }}\left(
z_{1},\frac{z_{2}}{a^{t_{2}-t_{1}}},\frac{z_{3}}{a^{t_{3}-t_{1}}}\right)
\times \left( 1-a^{J+t_{1}-t_{3}+1}\right) \right) .
\label{Eq_Distribution_Function_2}
\end{align}%
A similar calculation yields 
\begin{align}
& \mathds{P}\Bigg (\bigvee_{j=t_{1}-t_{2}}^{-1}\left\{
a^{j+t_{2}-t_{1}}(1-a)Z(t_{1}-j,\mathbf{x}_{2}-(j+t_{2}-t_{1})\boldsymbol{%
\tau })\right\} \leq z_{2},  \notag \\
& \ \ \ \ \  \bigvee_{j=t_{1}-t_{2}}^{-1}\left\{ a^{j+t_{3}-t_{1}}(1-a)Z(t_{1}-j,\mathbf{x}%
_{3}-(j+t_{3}-t_{1})\boldsymbol{\tau })\right\} \leq z_{3}\Bigg ) \nonumber
\label{Eq_Distribution_Function_3} \\
& =\exp \left( -\mathcal{V}_{\mathbf{x}_{2},\mathbf{x}_{3}-(t_{3}-t_{2})\boldsymbol{%
\tau }}\left( \frac{z_{2}}{a^{t_{2}-t_{1}}},\frac{z_{3}}{a^{t_{3}-t_{1}}}%
\right) \times (1-a^{t_{2}-t_{1}})\right) .  
\end{align}%
Furthermore, we have that
\begin{align}
& \mathds{P}\left( \bigvee_{j=t_{1}-t_{3}}^{t_{1}-t_{2}-1}\left\{
a^{j+t_{3}-t_{1}}(1-a)Z(t_{1}-j,\mathbf{x}_{3}-(j+t_{3}-t_{1})\boldsymbol{%
\tau })\right\} \leq z_{3}\right)   \notag \\
& =\prod_{j=t_{1}-t_{3}}^{t_{1}-t_{2}-1}\mathds{P}\left( Z(t_{1}-j,\mathbf{x}%
_{3}-(j+t_{3}-t_{1})\boldsymbol{\tau })\leq \frac{z_{3}}{%
a^{j}a^{t_{3}-t_{1}}(1-a)}\right)   \notag \\
& =\exp \left( -\frac{a^{t_{3}-t_{1}}(1-a)}{z_{3}}%
\sum_{t_{1}-t_{3}}^{t_{1}-t_{2}-1}a^{j}\right)   \notag \\
& =\exp \left( -\frac{1-a^{t_{3}-t_{2}}}{z_{3}}\right) 
\label{Eq_Distribution_Function_4}
\end{align}%
and similarly%
\begin{equation*}
\mathds{P}\left( \bigvee_{J+t_{1}-t_{2}+1}^{J}a^{j}(1-a)Z(t_{1}-j,\mathbf{x}%
_{1}-j\boldsymbol{\tau })\leq z_{1}\right) =\exp \left( -\frac{%
a^{J+t_{1}-t_{2}+1}(1-a^{t_{2}-t_{1}})}{z_{1}}\right) .
\end{equation*}%
Finally, a similar calculation as in \eqref{Eq_Distribution_Function_2}
yields 
\begin{align}
& \mathds{P}\Bigg( \bigvee_{j=J+t_{1}-t_{3}+1}^{J+t_{1}-t_{2}}\left\{
a^{j}(1-a)Z(t_{1}-j,\mathbf{x}_{1}-j\boldsymbol{\tau })\right\} \leq
z_{1}, \nonumber
\\& \ \ \ \ \ \bigvee_{j=J+t_{1}-t_{3}+1}^{J+t_{1}-t_{2}}\left\{
a^{j+t_{2}-t_{1}}(1-a)Z(t_{1}-j,\mathbf{x}_{2}-(j+t_{2}-t_{1})\boldsymbol{%
\tau })\right\} \leq z_{2} \Bigg)   \notag \\
& =\exp \left( -\mathcal{V}_{\mathbf{x}_{1},\mathbf{x}_{2}+(t_{2}-t_{1})%
\boldsymbol{\tau }}\left( z_{1},\frac{z_{2}}{a^{t_{2}-t_{1}}}\right)
a^{J+t_{1}-t_{3}+1}(1-a^{t_{3}-t_{2}})\right) .
\label{Eq_Distribution_Function_6}
\end{align}%
Inserting \eqref{Eq_Distribution_Function_2}, %
\eqref{Eq_Distribution_Function_3}, \eqref{Eq_Distribution_Function_4} and %
\eqref{Eq_Distribution_Function_6} in \eqref{Eq_Distribution_Function_1}, we
obtain, since $\lim_{J\rightarrow \infty }a^{J}=0$, 
\begin{align*}
& \mathds{P}(X(t_{1},\mathbf{x}_{1})\leq z_{1},X(t_{2},\mathbf{x}_{2})\leq
z_{2},X(t_{3},\mathbf{x}_{3})\leq z_{3}) \\
& =\lim_{J\rightarrow \infty }\mathds{P}\Bigg(\bigvee_{j=0}^{J}\left\{
a^{j}(1-a)Z(t_{1}-j,\mathbf{x}_{1}-j\boldsymbol{\tau })\leq \right\}
z_{1},\bigvee_{j=0}^{J}\left\{ a^{j}(1-a)Z(t_{2}-j,\mathbf{x}_{2}-j\boldsymbol{%
\tau })\right\} \leq z_{2}, \\
& \ \ \ \ \ \ \ \ \ \bigvee_{j=0}^{J}\left\{ a^{j}(1-a)Z(t_{3}-j,\mathbf{x}%
_{3}-j\boldsymbol{\tau })\right\} \leq z_{3}\Bigg) \\
& =\exp \left( - \mathcal{V}_{\mathbf{x}_{1},\mathbf{x}_{2}-(t_{2}-t_{1})%
\boldsymbol{\tau },\mathbf{x}_{3}-(t_{3}-t_{1})\boldsymbol{\tau }}\left(
z_{1},\frac{z_{2}}{a^{t_{2}-t_{1}}},\frac{z_{3}}{a^{t_{3}-t_{1}}}\right) \right)  \\
& \ \ \ \times \exp \left( -\mathcal{V}_{\mathbf{x}_{2},\mathbf{x}%
_{3}-(t_{3}-t_{2})\boldsymbol{\tau }}\left( \frac{z_{2}}{a^{t_{2}-t_{1}}},%
\frac{z_{3}}{a^{t_{3}-t_{1}}}\right) \times (1-a^{t_{2}-t_{1}})\right)
\times \exp \left( -\frac{1-a^{t_{3}-t_{2}}}{z_{3}}\right) .
\end{align*}
\end{proof}

\subsection{For Proposition \protect\ref%
{Prop_Spatio_Temporal_Extremal_Coefficient}}

\begin{proof}
Applying \eqref{Eq_General_Distribution_Function} with $M=2$ and setting $%
z_{1}=z_{2}=u$ for $u>0$, we obtain 
\begin{align*}
\mathds{P}(X(t_{1},\mathbf{x}_{1}))\leq u,X(t_{2},\mathbf{x}_{2}) \leq u)&
=\exp \left( -\mathcal{V}_{\mathbf{x}_{1},\mathbf{x}_{2}-(t_{2}-t_{1})%
\boldsymbol{\tau }}\left( u,\frac{u}{a^{t_{2}-t_{1}}}\right) \right) \exp
\left( -\frac{1-a^{t_{2}-t_{1}}}{u}\right)  \\
& =\exp \left( -\frac{\mathcal{V}_{\mathbf{x}_{1},\mathbf{x}%
_{2}-(t_{2}-t_{1})\boldsymbol{\tau }}\left( 1,a^{t_{1}-t_{2}}\right)
+1-a^{t_{2}-t_{1}}}{u}\right),
\end{align*}
yielding \eqref{Eq_Spatio_Temporal_Extremal_Coefficient} by definition of the spatio-temporal extremal coefficient. \\
In the same way as in the purely spatial case \citep[see e.g.][p.379]{cooley2006variograms}, it is easy to show the following link between the spatio-temporal $\Phi_1$-madogram and the spatio-temporal extremal coefficient:
\begin{equation}
\label{Eq_F_Madogram_Bis}
\nu_F(t_2-t_1, \Mb{x}_2-\Mb{x}_1)=\frac{1}{2} \frac{\theta(t_2-t_1, \Mb{x}_2-\Mb{x}_1)-1}{\theta(t_2-t_1, \Mb{x}_2-\Mb{x}_1)+1}=\frac{1}{2}-\frac{1}{\theta(t_2-t_1, \Mb{x}_2-\Mb{x}_1)+1}.
\end{equation}
Inserting \eqref{Eq_Spatio_Temporal_Extremal_Coefficient} into \eqref{Eq_F_Madogram_Bis} gives \eqref{Eq_F_Madogram}.
\end{proof}

\subsection{For Theorem \protect\ref{Prop_Geometric_Ergodicity}}

\begin{proof}
We show that the two assumptions appearing in Theorem 3.6 in \cite{hairer2010p} are satisfied. 

\medskip 

\noindent \textbf{Assumption 1.}
\textit{There exists a function $L:\mathcal{C}_{%
\mathds{S}^{2}}\rightarrow \lbrack 0,\infty )$ and constants $K\geq 0$ and $%
\gamma \in (0,1)$ such that%
\begin{equation*}
\mathds{E}[L(X(t,\mathbf{\cdot }))|X(t-1,\mathbf{\cdot })=h(\cdot )]\leq
\gamma L(h)+K
\end{equation*}%
for all $h\in \mathcal{C}_{\mathds{S}^{2}}$.
} 

\medskip

\noindent Let us choose $L(h)=\left\Vert h^{\gamma
}\right\Vert _{\infty }$ with $0<\gamma <1$. 
\\ First, note that the $\left( U_{i}\right)_{i\geq 1}$ in \eqref{DistributionZspatial} satisfy $\left( U_{i}\right)_{i\geq 1}=\left( P_{i}^{-1}\right) _{i\geq 1}$,
where $\left( P_{i}\right) _{i\geq 1}$ are the points of an homogeneous
Poisson point process on $\mathds{R}_{+}$ with constant intensity equal to
one. Hence, the highest $U_i$ corresponds to the smallest $P_i$, which is exponentially distributed with parameter 1. Hence, its inverse follows the standard Fr\'echet distribution. Moreover, for $f$ defined in \eqref{Eq_von_Mises_Fisher}, $\| f \|_{\infty}$ is reached for $\Mb{x}=\boldsymbol{\mu}$ and is finite. Therefore, the process $Z$ defined in \eqref{DistributionZspatial} satisfies
\begin{equation}
\label{Eq_1_Proof_Th4}
\left\Vert Z^{\gamma }\right\Vert _{\infty }\overset{d}{=}Y^{\gamma
}\left\Vert f^{\gamma }\right\Vert _{\infty },
\end{equation}%
where $Y$ is a random variable with standard Fr\'{e}chet distribution. Moreover,
\begin{eqnarray}
\left\Vert \max (ah(R_{\theta, \Mb{u} }\mathbf{x}),(1-a)Z(t, \mathbf{x}))^{\gamma
}\right\Vert _{\infty } &=&\left\Vert \max (a^{\gamma }h^{\gamma }(R_{\theta, \Mb{u}}\mathbf{x}),(1-a)^{\gamma }Z^{\gamma }(t, \mathbf{x}))\right\Vert _{\infty } \nonumber \\
&\leq &\max (a^{\gamma }\left\Vert h^{\gamma }\right\Vert _{\infty
},(1-a)^{\gamma }\left\Vert Z^{\gamma }\right\Vert _{\infty }).
\label{Eq_2_Proof_Th4}
\end{eqnarray}
Using \eqref{Eq_1_Proof_Th4} and \eqref{Eq_2_Proof_Th4}, we obtain
\begin{align*}
&\mathds{E}\left[ L(X(t,\mathbf{\cdot }))|X(t-1,\mathbf{\cdot })=h(\cdot )%
\right]  
\\&=\mathds{E}[\left\Vert \max (a^{\gamma }h^{\gamma }(R_{\theta, \Mb{u} }\mathbf{x}%
),(1-a)^{\gamma }Z^{\gamma }(t, \mathbf{x})) \right\Vert _{\infty }] 
\\&\leq \mathds{E} [\max (a^{\gamma }\left\Vert h^{\gamma }\right\Vert
_{\infty },(1-a)^{\gamma }\left\Vert Z^{\gamma }\right\Vert _{\infty }) ] 
\\&= a^{\gamma }\left\Vert h^{\gamma }\right\Vert _{\infty }\mathds{P} \left(
Y^{\gamma } \leq \frac{ a^{\gamma }\left\Vert h^{\gamma }\right\Vert _{\infty
}} {(1-a)^{\gamma }\left\Vert f^{\gamma } \right\Vert_{\infty} } \right) +(1-a)^{\gamma } \| f\|_{\infty}
\mathds{E}  \left[ Y^{\gamma }\mathds{I}_{ \left \{ Y^{\gamma }\geq \frac{ a^{\gamma }\left\Vert
h^{\gamma }\right\Vert _{\infty }} { (1-a)^{\gamma }\left\Vert f^{\gamma
}\right\Vert_{\infty} } \right \} } \right] 
\\&\leq a^{\gamma }\left\Vert h^{\gamma }\right\Vert _{\infty }+(1-a)^{\gamma
} \| f\|_{\infty} \Gamma (1-\gamma ) 
\\ &=a^{\gamma }L(h)+(1-a)^{\gamma } \| f\|_{\infty} \Gamma (1-\gamma ),
\end{align*}
yielding Assumption 1.

\medskip 

\noindent \textbf{Assumption 2.}
\textit{
We denote by $\| . \|_{TV}$ the total variation distance between two probability measures.
For every $R>0$, there exists a constant $%
\alpha >0$ such that
\begin{equation*}
\sup_{h_{1},h_{2}\in D}\left\Vert \mathcal{P}_{h_{1}}-\mathcal{P}%
_{h_{2}}\right\Vert _{TV}\leq 2(1-\alpha ),
\end{equation*}%
where $D=\{h_{1},h_{2}:L(h_{1})+L(h_{2})\leq R\}$, or equivalently%
\begin{equation*}
\sup_{h_{1},h_{2}\in D,\left\Vert \varphi \right\Vert _{\infty }\leq
1}\left\vert \mathds{E}\left[ \varphi (X(t,\mathbf{\cdot }))|X(t-1,\mathbf{%
\cdot })=h_{1}(\cdot )\right] -\mathds{E}\left[ \varphi (X(t,\mathbf{\cdot }%
))|X(t-1,\mathbf{\cdot })=h_{2}(\cdot )\right] \right\vert \leq 2(1-\alpha ).
\end{equation*}
}

\medskip

\noindent Using \eqref{Eq_Recurrence_Model_Type_2}, we have that
\begin{align}
&\left\vert \mathds{E}\left[ \varphi (X(t,\mathbf{\cdot }))|X(t-1,\mathbf{%
\cdot })=h_{1}(\cdot )\right] -\mathds{E}\left[ \varphi (X(t,\mathbf{\cdot }%
))|X(t-1,\mathbf{\cdot })=h_{2}(\cdot )\right] \right\vert  \nonumber
\\&= \left\vert \mathds{E}\left[ \varphi (\max ( ah_{1}(R_{\theta, \Mb{u} }\mathbf{x}%
),(1-a)Z(t, \mathbf{x})))-\varphi (\max (ah_{2}(R_{\theta, \Mb{u} }\mathbf{x}),(1-a)Z(t, 
\mathbf{x})))\right] \right\vert  \nonumber
\\&\leq \mathds{E}\left[ \left\vert \varphi (\max (ah_{1}(R_{\theta, \Mb{u} }\mathbf{x%
}),(1-a)Z(t, \mathbf{x})))-\varphi (\max (ah_{2}(R_{\theta, \Mb{u} }\mathbf{x}),(1-a)Z(t,
\mathbf{x})))\right\vert \right].
\label{Eq_3_Proof_Th4}
\end{align}
Moreover, for $b,c,z\in \mathcal{C}_{\mathds{S}^{2}}$, we know that
\begin{align*}
&\left\vert \varphi (\max (b(\Mb{x}),z(\Mb{x}) ))-\varphi (\max (c(\Mb{x}),z(\Mb{x})))\right\vert =0 \ \forall \Mb{x} \in \mathds{S}^2
\\& \mbox{ if and only if } \max( b(\mathbf{x}), c(\mathbf{x})) \leq z(\mathbf{x}) \ \forall \mathbf{x}%
\in \mathds{S}^{2}. 
\end{align*}
Therefore, for functions $\varphi$ satisfying $\| \varphi \|_{\infty} \leq 1$, 
\begin{align}
&\mathds{E}\left[ \left\vert \varphi (\max (ah_{1}(R_{\theta, \Mb{u} }\mathbf{x}%
),(1-a)Z(t, \mathbf{x})))-\varphi (\max (ah_{2}(R_{\theta, \Mb{u} }\mathbf{x}),(1-a)Z(t,
\mathbf{x})))\right\vert \right] \nonumber  
\\&\leq 2\left( 1-\mathds{P}\left( \bigcap_{\Mb{x} \in \mathds{S}^{2}} \left \{ Z(t, \mathbf{x})\geq \frac{a}{%
(1-a)} \max( h_{1}(R_{\theta, \Mb{u} }\mathbf{x}), h_{2}(R_{\theta, \Mb{u} }\mathbf{x}) ) \right \} \right)
\right)  \nonumber
\\&= 2\left( 1-\mathds{P}\left( \bigcap_{\Mb{x} \in \mathds{S}^{2}} \left \{ Z(t, R_{-\theta, \Mb{u} }\mathbf{x})\geq
\frac{a}{(1-a)} \max( h_{1}(\mathbf{x}), h_{2}(\mathbf{x})) \right \} \right) \right)  \nonumber
\\&=2\left( 1-\mathds{P}\left( \bigcap_{\Mb{x} \in \mathds{S}^{2}} \left \{ Z(t, \mathbf{x})\geq \frac{a}{(1-a)}%
\max( h_{1}(\mathbf{x}), h_{2}(\mathbf{x})) \right \} \right) \right) \nonumber
\\& \leq 2 \left(1 - \mathds{P}\left( \bigwedge_{\mathbf{%
x} \in \mathds{S}^{2}} \{ Z(t, \mathbf{x}) \} \geq \frac{a}{(1-a)} \max( \Vert h_{1} \Vert _{\infty
}, \left\Vert h_{2}\right\Vert _{\infty }) \right) \right).
\label{Eq_4_Proof_Th4}
\end{align}%
The quantity $\bigwedge_{\mathbf{x}\in \mathds{S}^{2}} \{ f(\Mb{x}; \boldsymbol{\mu}_1, \kappa) \}$ is reached for $\Mb{x}=-\boldsymbol{\mu}_1$. Thus, 
\begin{align*}
Z(\mathbf{x})=\bigvee_{i=1}^{\infty } \left \{ U_{i}f(\mathbf{x};\boldsymbol{\mu }%
_{i},\kappa ) \right \} \geq P_{1}^{-1}f(\mathbf{x};\mathbf{\mu }_{1},\kappa ) &\geq
P_{1}^{-1} \bigwedge_{\mathbf{x}\in \mathds{S}^{2}} \{ f(\mathbf{x};\boldsymbol{\mu }%
_{1},\kappa ) \}
\\& =P_{1}^{-1}\frac{\kappa }{4\pi \sinh \kappa }\exp \left(
-\kappa \right).
\end{align*}%
It follows that%
\begin{align}
&\mathds{P}\left( \bigwedge_{\mathbf{x} \in \mathds{S}^{2}} \{ Z(\mathbf{x}) \} \geq \frac{a}{(1-a)}%
\max( \left \Vert h_{1}\right\Vert _{\infty }, \left\Vert h_{2}\right\Vert
_{\infty } ) \right) \nonumber
\\&\geq \mathds{P}\left( P_{1}^{-1}\frac{\kappa }{4\pi \sinh \kappa }\exp
\left( -\kappa \right) \geq \frac{a}{(1-a)} \max( \left\Vert h_{1}\right\Vert
_{\infty }, \left\Vert h_{2}\right\Vert _{\infty }) \right)  \nonumber
\\& \geq \mathds{P}\left( P_{1}^{-1}\geq \frac{4\pi \sinh \kappa a}{\kappa
(1-a)}\exp \left( \kappa \right) R\right),
\label{Eq_5_Proof_Th4}
\end{align}%
noting that $\max( \left\Vert h_{1}\right\Vert _{\infty }, \left\Vert
h_{2}\right\Vert _{\infty } ) \leq R$. \\
Therefore, combining \eqref{Eq_3_Proof_Th4}, \eqref{Eq_4_Proof_Th4} and \eqref{Eq_5_Proof_Th4}, we obtain
\begin{align*}
& \sup_{h_{1},h_{2}\in D,\left\Vert \varphi \right\Vert _{\infty }\leq
1}\left\vert \mathds{E}\left[ \varphi (X(t,\mathbf{\cdot }))|X(t-1,\mathbf{%
\cdot })=h_{1}(\cdot )\right] -\mathds{E}\left[ \varphi (X(t,\mathbf{\cdot }%
))|X(t-1,\mathbf{\cdot })=h_{2}(\cdot )\right] \right\vert
\\& \leq 2 \left( 1- \mathds{P}\left( P_{1}^{-1}\geq \frac{4\pi \sinh \kappa a}{\kappa
(1-a)}\exp \left( \kappa \right) R\right) \right)=2(1-\alpha),
\end{align*}
denoting $$\alpha=\mathds{P}\left( P_{1}^{-1}\geq \frac{4\pi \sinh \kappa a}{\kappa (1-a)}\exp \left( \kappa \right) R\right)>0.$$
Hence, Assumption 2 holds. 

\medskip

\noindent Finally, the application of Theorem 3.6 in \cite{hairer2010p} yields the result.
\end{proof}

\subsection{For Proposition \protect\ref{Prop_Bivariate_Density_MARMA_01}}

\begin{proof}
We have that 
\begin{align*}
& f_{(t_{1},\mathbf{x}_{1}),(t_{2},\mathbf{x}_{2})}(z_{1},z_{2},\boldsymbol{%
\theta }) \\
& =\frac{\partial }{\partial z_{1}}\frac{\partial }{\partial z_{2}}\exp
\left( -\mathcal{V}_{\mathbf{x}_{1},\mathbf{x}_{2}-(t_{2}-t_{1})\boldsymbol{%
\tau }}\left( z_{1},\frac{z_{2}}{a^{t_{2}-t_{1}}}\right) -\frac{1-\phi
^{t_{2}-t_{1}}}{z_{2}}\right)  \\
& =\frac{\partial }{\partial z_{1}}\Bigg( \exp \left( -\mathcal{V}_{\mathbf{x%
}_{1},\mathbf{x}_{2}-(t_{2}-t_{1})\boldsymbol{\tau }}\left( z_{1},\frac{z_{2}%
}{a^{t_{2}-t_{1}}}\right) -\frac{1-\phi ^{t_{2}-t_{1}}}{z_{2}}\right) \times
\Big( -\frac{\partial }{\partial z_{2}}\mathcal{V}_{\mathbf{x}_{1},\mathbf{x%
}_{2}-(t_{2}-t_{1})\boldsymbol{\tau }}\left( z_{1},\frac{z_{2}}{%
a^{t_{2}-t_{1}}}\right) 
\\& \ \ \ +\frac{1-a^{t_{2}-t_{1}}}{z_{2}^{2}}\Big) \Bigg) 
\\
& =\exp \left( -\mathcal{V}_{\mathbf{x}_{1},\mathbf{x}_{2}-(t_{2}-t_{1})%
\boldsymbol{\tau }}\left( z_{1},\frac{z_{2}}{a^{t_{2}-t_{1}}}\right) -\frac{%
1-\phi ^{t_{2}-t_{1}}}{z_{2}}\right) \times \left( -\frac{\partial }{%
\partial z_{1}}\mathcal{V}_{\mathbf{x}_{1},\mathbf{x}_{2}-(t_{2}-t_{1})%
\boldsymbol{\tau }}\left( z_{1},\frac{z_{2}}{a^{t_{2}-t_{1}}}\right) \right) 
\\
& \ \ \ \times \left( -\frac{\partial }{\partial z_{2}}\mathcal{V}_{\mathbf{x%
}_{1},\mathbf{x}_{2}-(t_{2}-t_{1})\boldsymbol{\tau }}\left( z_{1},\frac{z_{2}%
}{a^{t_{2}-t_{1}}}\right) +\frac{1-a^{t_{2}-t_{1}}}{z_{2}^{2}}\right)  \\
& \ \ \ +\exp \left( -\mathcal{V}_{\mathbf{x}_{1},\mathbf{x}%
_{2}-(t_{2}-t_{1})\boldsymbol{\tau }}\left( z_{1},\frac{z_{2}}{%
a^{t_{2}-t_{1}}}\right) -\frac{1-\phi ^{t_{2}-t_{1}}}{z_{2}}\right) \times
\left( -\frac{\partial }{\partial z_{1}}\frac{\partial }{\partial z_{2}}%
\mathcal{V}_{\mathbf{x}_{1},\mathbf{x}_{2}-(t_{2}-t_{1})\boldsymbol{\tau }%
}\left( z_{1},\frac{z_{2}}{a^{t_{2}-t_{1}}}\right) \right) ,
\end{align*}%
yielding the result.
\end{proof}

\subsection{For Corollary \protect\ref%
{Prop_Bivariate_Density_MARMA_01_Smith}}

\begin{proof}
We denote by $w=\frac{h}{2}+\frac{1}{h}\log \left( \frac{z_{2}}{z_{1}}%
\right) $ and $v=h-w$, where $h=\sqrt{(\mathbf{x}_{2}-\mathbf{x}%
_{1})^{\prime } \Sigma^{-1}(\mathbf{x}_{2}-\mathbf{x}_{1})}$. From 
\cite{padoan2010likelihood}, p.275, we know that 
\begin{equation*}
-\frac{\partial }{\partial z_{1}}\mathcal{V}_{\mathbf{x}_{1},\mathbf{x}%
_{2}}(z_{1},z_{2})=\frac{\Phi (w)}{z_{1}^{2}}+\frac{\phi (w)}{hz_{1}^{2}}-%
\frac{\phi (v)}{hz_{1}z_{2}},\quad -\frac{\partial }{\partial z_{2}}\mathcal{%
V}_{\mathbf{x}_{1},\mathbf{x}_{2}}(z_{1},z_{2})=\frac{\Phi (v)}{z_{2}^{2}}+%
\frac{\phi (v)}{hz_{2}^{2}}-\frac{\phi (w)}{hz_{1}z_{2}}\mbox{ and }
\end{equation*}%
\begin{equation*}
-\frac{\partial ^{2}}{\partial z_{1}\partial z_{2}}\mathcal{V}_{\mathbf{x}%
_{1},\mathbf{x}_{2}}(z_{1},z_{2})=\frac{v\phi (w)}{h^{2}z_{1}^{2}z_{2}}+%
\frac{w\phi (v)}{h^{2}z_{1}z_{2}^{2}}.
\end{equation*}%
Hence, we obtain 
\begin{equation}
-\frac{\partial }{\partial z_{1}}\mathcal{V}_{\mathbf{x}_{1},\mathbf{x}%
_{2}-(t_{2}-t_{1})\boldsymbol{\tau }}\left( z_{1},\frac{z_{2}}{%
a^{t_{2}-t_{1}}}\right) =\frac{\Phi (w_{1})}{z_{1}^{2}}+\frac{\phi (w_{1})}{%
h_{1}z_{1}^{2}}-\frac{a^{t_{2}-t_{1}}\phi (v_{1})}{h_{1}z_{1}z_{2}},
\label{Bivariate_Density_Eq_1}
\end{equation}%
\begin{align}
-\frac{\partial }{\partial z_{2}}\mathcal{V}_{\mathbf{x}_{1},\mathbf{x}%
_{2}-(t_{2}-t_{1})\boldsymbol{\tau }}\left( z_{1},\frac{z_{2}}{%
a^{t_{2}-t_{1}}}\right) & =\frac{1}{a^{t_{2}-t_{1}}}\left( \frac{%
a^{2(t_{2}-t_{1})}\Phi (v_{1})}{z_{2}^{2}}+\frac{a^{2(t_{2}-t_{1})}\phi
(v_{1})}{h_{1}z_{2}^{2}}-\frac{a^{t_{2}-t_{1}}\phi (w_{1})}{h_{1}z_{1}z_{2}}%
\right)   \notag \\
& =\frac{a^{t_{2}-t_{1}}\Phi (v_{1})}{z_{2}^{2}}+\frac{a^{t_{2}-t_{1}}\phi
(v_{1})}{h_{1}z_{2}^{2}}-\frac{\phi (w_{1})}{h_{1}z_{1}z_{2}}
\label{Bivariate_Density_Eq_2}
\end{align}%
and 
\begin{align}
-\frac{\partial ^{2}}{\partial z_{1}\partial z_{2}}\mathcal{V}_{\mathbf{x}%
_{1},\mathbf{x}_{2}-(t_{2}-t_{1})\boldsymbol{\tau }}\left( z_{1},\frac{z_{2}%
}{a^{t_{2}-t_{1}}}\right) &=\frac{1}{a^{t_{2}-t_{1}}}\left( \frac{%
a^{t_{2}-t_{1}}v_{1}\phi (w_{1})}{h_{1}^{2}z_{1}^{2}z_{2}}+\frac{%
a^{2(t_{2}-t_{1})}w_{1}\phi (v_{1})}{h_{1}^{2}z_{1}z_{2}^{2}}\right) \nonumber
\\& =\frac{%
v_{1}\phi (w_{1})}{h_{1}^{2}z_{1}^{2}z_{2}}+\frac{a^{t_{2}-t_{1}}w_{1}\phi
(v_{1})}{h_{1}^{2}z_{1}z_{2}^{2}}.  \label{Bivariate_Density_Eq_3}
\end{align}%
Finally, we have that 
\begin{equation}
\mathcal{V}_{\mathbf{x}_{1},\mathbf{x}_{2}-(t_{2}-t_{1})\boldsymbol{\tau }%
}\left( z_{1},\frac{z_{2}}{a^{t_{2}-t_{1}}}\right) =\frac{\Phi (w_{1})}{z_{1}%
}+\frac{a^{t_{2}-t_{1}}\Phi (v_{1})}{z_{2}}.  \label{Bivariate_Density_Eq_4}
\end{equation}%
Inserting \eqref{Bivariate_Density_Eq_1}, \eqref{Bivariate_Density_Eq_2}, %
\eqref{Bivariate_Density_Eq_3} and \eqref{Bivariate_Density_Eq_4} in %
\eqref{Eq_Bivariate_Density_MARMA_01}, we obtain the result.
\end{proof}

\newpage
\bibliographystyle{apalike}
\bibliography{Extension_MARMA_Spatial}

\end{document}